\documentclass[a4]{article}
\usepackage{fullpage}
\usepackage{amsfonts}
\usepackage{amsmath}
\usepackage{graphicx}
\usepackage[mathscr]{eucal}

\title{On transversal connecting orbits of Lagrangian systems in non-stationary force field: Newton-Kantorovich approach}
\author{Alexey Ivanov}
\date{}
\begin{document}
\renewcommand{\theequation}{\arabic{section}.\arabic{equation}}
\maketitle

\begin{abstract}
We consider a natural Lagrangian system defined on a complete Riemannian 
manifold being subjected to the action of a non-stationary force field with potential $U(q,t) = f(t)V(q)$. It is assumed that the factor $f(t)$ tends to $\infty$ as $t\to \pm\infty$ and vanishes at a unique point $t_{0}\in \mathbb{R}$. Let $X_{+}$, $X_{-}$ denote the sets of isolated critical points of $V(x)$ at which $U(x,t)$ as a function of $x$ attains its maximum for any fixed $t> t_{0}$ and $t<t_{0}$, respectively. Under nondegeneracy conditions on points of $X_{\pm}$ we apply Newton-Kantorovich type method to study the existence of transversal doubly asymptotic trajectories connecting $X_{-}$ and $X_{+}$. Conditions on the Riemannian manifold and the potential which guarantee the existence of such orbits are presented. Such connecting trajectories are obatained by continuation of geodesics defined in a vicinity of the point $t_{0}$ to the whole real line.
\end{abstract}



\section{Introduction}

Since pioneer works of Poincar\'e the doubly asymptotic  (i.e. homoclinic and heteroclinic) trajectories became a subject of many scientific papers. Such trajectories belong to the intersection of invariant (stable and unstable) manifolds associated to hyperbolic objects (equilibria, periodic orbits, invariant tori). It was found \cite{Poi} that the mutual disposition of the invariant manifolds with homoclinic or heteroclinic intersections in the phase space of a dynamical system can be extremely complicated. This leads to a complex behaviour of the system in a vicinity of such invariant manifolds. In particular, it was proved by Birkhoff \cite{Bir} and Smale \cite{Sma} that existence of a transversal connecting orbit implies chaotic dynamics of a system. 

Different methods and techniques are used in studying homoclinic and heteroclinic trajectories. They can be divided into three main parts: asymptotic (or perturbative), variational and numeric ones. Asymptotic methods (Melnikov's method \cite{Mel}, exponentially small splitting methods \cite{GelLaz}, \cite{Tre}, singular perturbation methods \cite{Fen}, \cite{Wig} and others) allow establishing the existence of such trajectories via analysis of invariant manifolds of hyperbolic objects for the unperturbed system. These methods resolve not only the existence problem, but also provide  more information on the geometry of invariant manifolds and its intersections for the perturbed system. However, knowledge of the unperturbed system is crucial for such kind of methods. In contrast, variational methods \cite{BolKoz}, \cite{CotiRab}, \cite{Rab}, \cite{Ang} can be applied in much more general context. One of the consequences of such generality is the absence of information on the transversality of the obtained connecting orbits. Methods of the third kind takes some intermediate position. Using different techniques (Newton's method, shadowing \cite{Pal}, \cite{CKP} etc.) they can be applied for sufficiently general systems to construct doubly asymptotic trajectories and check their transversality. But these methods as asymptotic ones need some {\it a priori} information: initial approximation for the Newton's method or pseudo-orbits for shadowing. To obtain such initial data one often uses numeric simulations.

In the present work we consider a time-dependent Lagrangian system defined on a complete Riemannian manifold. Particularly, let a compact Riemannian manifold $\mathcal{M}$ be the configuration space of a Lagrangian system with Lagrangian 
$L\in C^{3}(T\mathcal{M} \times {\mathbb{R}}, {\mathbb{R}})$ such that 
\begin{eqnarray}
L(q,\dot q, t)  = K(q,\dot q) - U(q,t).
\end{eqnarray}
It is assumed that the kinetic energy $K\in C^{3}(T\mathcal{M}, {\mathbb{R}})$ is
a positive definite quadratic form in velocity $\dot q$ and the potential 
$U(q,t)\in C^{3}(\mathcal{M} \times {\mathbb{R}}, {\mathbb{R}})$ has the following representation 
\begin{eqnarray}
U(q,t) = f(t)V(q).
\end{eqnarray}
In addition it is supposed that the factor $f$ satisfies the following assumptions:
\newline
\newline
($A_{1}$) there exists a unique $t_{0}\in \mathbb{R}$ such that $f(t_{0}) = 0$,
\newline
\newline
($A_{2}$) $\vert f(t)\vert \to +\infty$ as $t\to \pm\infty$,
\newline
\newline
($A_{3}$) $f''(t) f(t) < 3 (f'(t))^{2}/2$ for all $t\in \mathbb{R}$.
\newline
\newline
Without lost of generality one may always suppose that $t_{0} = 0$. Hence, the factor $f(t)$ does not change the sign on the intervals $\mathbb{R}_{\pm}$, where $\mathbb{R}_{+} = (0, +\infty)$ and  $\mathbb{R}_{-} = (-\infty, 0)$.

Such systems arise in different areas of physics \cite{Koz91}, \cite{BertBol}, \cite{BorKozMam}. One of the applications motivating the present work lies in the theory of Lagrangian systems with turning points. If a Lagrangian system has a turning point of order $m$ at $t=t_{0}$ then in a small vicinity of $t_{0}$ it may be approximated \cite{Iva17} by the so-called reference system of the type (1.1) with potential $U(q,t) = f(t)V(q)$ and the factor $f(t) = (t-t_{0})^{m}$. It was proved \cite{Iva17} that if additionally such Lagrangian system with turning points is time-periodic then in the adiabatic limit, i.e. when the Lagrangian is of the form $L = L(q, \dot q, \varepsilon t)$ and the parameter $\varepsilon\ll 1$, it possesses  plenty of connecting orbits and multi-bump trajectories. However, the proof of this result is based on an assumption that the reference system associated to the turning point $t_{0}$ has infinitely many transversal connecting trajectories. Using variational arguments one may prove \cite{Iva16} that a system of type (1.1) indeed possesses infinitley many connecting orbits, but the transversality assumption still has to be checked. The aim of this paper is to obtain sufficient conditions on the manifold $\mathcal{M}$ and the potential $U$ which guarantee the existence of transversal connecting orbits of the system (1.1). One may note that the system (1.1) is parameter free, thus, neither asymptotic nor variational methods cannot be applied for this purpose. 

Since $\mathcal{M}$ is a compact manifold, the function $V$ has minimum and maximum on  $\mathcal{M}$. Following notations of \cite{Iva16} for any fixed $t>0$ (resp. $t<0$) we denote by $X_{+}(t)$ (resp. $X_{-}(t)$) the subset of $\mathcal{M}$ on which the potential $U(x,t)$ considered as a function of the variable $x$ attains its maximum. The condition $A_{1}$ implies that the factor $f$ may change sign only at $t=0$. Taking this into account and using the representation (1.2) we may conclude that the subset $X_{+}(t)$ (resp. $X_{-}(t)$) does not depend on $t$ on the interval  $\mathbb{R}_{+}$  (resp. $\mathbb{R}_{-}$). Hence, one may skip the dependence on $t$ in the definition of the subsets $X_{\pm}$.  We suppose that \newline
\newline
($A_{4}$) $X_{\pm}$ consist of nondegenerate isolated critical points of $V$.
\newline

We will say that a solution $q:\mathbb{R}\to \mathcal{M}$ is a {\it heteroclinic} ({\it homoclinic})  solution if there exist 
$x_{-}, x_{+}\in \mathcal{M}$ (for homoclinic solution $x_{-}=x_{+}$) such that $q$ joins $x_{-}$ to $x_{+}$, i.e. $\lim\limits_{t\to \pm\infty}q(t) = x_{\pm}$ and $\lim\limits_{t\to \pm\infty}\dot q(t) = 0$. 

The main pecularity of the system (1.1) is vanishing of the potential $U(q,t)$ at $t=0$. If we fix two points $x_{\pm}\in X_{\pm}$ then for sufficiently small $T_{0}>0$ a solution 
$q_{0}: [-T_{0}, T_{0}]\to \mathcal{M}$ such that $q_{0}(\pm T_{0}) = x_{\pm}$ can be approximated by a geodesic which connects $x_{-}$ and $x_{+}$. Using Newton-Kantorovich method one may try to construct the solution $q_{0}$ and check its transversality. Then one may consider some $T_{1} > T_{0}$ and try to construct a solution $q_{1}: [-T_{1}, T_{1}]\to \mathcal{M}, q_{1}(\pm T_{1}) = x_{\pm}$ using $q_{0}$ as initial approximation. Continuing in the same way we obtain sequences of expanding intervals $[-T_{k},T_{k}]$ and solutions $q_{k}$. Since the factor $f$ tends to infinity as $\vert t \vert\to +\infty$ one may hope that the period of time which $q_{k}$ spends in a small (but fixed) neighborhood of the points $x_{\pm}$ increases with $k$. This will lead to fast convergence of the defined procedure, namely, $T_{k}\to +\infty$ and $q_{k}\to q_{\infty}$ as $k\to +\infty$, where $q_{\infty}$ is a transversal connecting orbit joining $x_{-}$ and $x_{+}$. 

The paper is organized as follows. In section 2 we introduce a Hilbert manifold of curves and define on this manifold the action functional whose critical points correspond to the doubly asymptotic trajectories of the system (1.1). Section 3 establishes expressions for the first and the second derivatives of the action functional and describes a relations between transversality of the connecting orbits and non-degeneracy of critical points. Section 4 is devoted to the Newton-Kantorovich theorem and its application to the present setting. In section 5 we describe a procedure for constructing transversal doubly asymptotic trajectories and provide conditions for its convergence. Finally, in section 6 we study a special case of the potential $U$, which  corresponds to the factor $f(t)=t^{m}, m\in \mathbb{N}$.

\section{Hilbert manifold, action functional, critical points}

\setcounter{equation}{0}

Consider a smooth embedding of the manifold $\mathcal{M}$ into $\mathbb{R}^{N}$ for $N=2n+1$ with $n = \dim \mathcal{M}$ and denote by $\langle \cdot, \cdot \rangle$ the Euclidean structure in $\mathbb{R}^{N}$ together with its restriction to $\mathcal{M}$. Let "$grad$" stands for the gradient operator with respect to the variable $x$. We fix two points $x_{\pm}\in X_{\pm}$ and will seek trajectories connecting $x_{-}$ and $x_{+}$. Note that in the case $X_{+} = X_{-}$ the points $x_{\pm}$ may coincide.

To simplify exposition we introduce a new time variable
\begin{equation*}
\xi = \int\limits_{0}^{t}\vert f(s)\vert^{1/2}{\rm d}s
\end{equation*}
and functions $r(\xi)$, $\sigma(\xi)$, $p(\xi)$:
\begin{equation*}
r(\xi) = \vert f(t)\vert^{1/2}\biggl\vert_{t=t(\xi)},\quad
\sigma(\xi) = {\rm sign}(f(t))\biggl\vert_{t=t(\xi)}, \quad
p(\xi) = \frac{r'(\xi)}{2 r(\xi)},
\end{equation*}
where $'$ denotes the derivative with respect to $\xi$.

One may note that $r$ is a non-negative $C^{1}$-function which has unique zero at the origin and satisfies the condition $r(\xi)\to \infty$ as $\vert \xi\vert \to \infty$. The function $p$ is defined on $\mathbb{R}\setminus \{0\}$ and due to assumption $A_{3}$ it decreases on $\mathbb{R}_{+}$ and increases on $\mathbb{R}_{-}$.The variable $t$ and the factor $f$ can be reconstructed via $r$ and $\sigma$ as 
$$
t = \int\limits_{0}^{\xi}\frac{{\rm d} s}{r(s)},\quad 
f(t) = \sigma(\xi) r^{2}(\xi)\biggl\vert_{\xi = \xi(t)}.
$$

We consider a vector space
\begin{equation*}
\mathcal{E}_{r} = \biggl\{v\in AC(\mathbb{R},\mathbb{R}^{N}): 
\|v\|_{r}^{2} = \int\limits_{\mathbb{R}} \bigl( |v'(\xi)|^{2} +  |v(\xi)|^{2}\bigr)r(\xi)
{\rm d}\xi < \infty\biggr\},
\end{equation*}
where $AC(\mathbb{R},\mathbb{R}^{N})$ is the set of absolutely continuous curves from $\mathbb{R}$ to $\mathbb{R}^{N}$.

Let $\mathcal{E}_{1} = W^{1,2}(\mathbb{R},\mathbb{R}^{N})$ be the Sobolev space with the norm $\| \cdot \|_{1}$ such that $\| v \|_{1}^{2} = \int\limits_{\mathbb{R}} \bigl( |v'(\xi)|^{2} + |v(\xi)|^{2}\bigr){\rm d}\xi$.

We summirize some results obtained in \cite{Iva16}. One may prove (lemma 1, \cite{Iva16}) that $\mathcal{E}_{r} \subset \mathcal{E}_{1}$  and 

\begin{eqnarray}
\| v \|_{1} \le C_{r} \| v \|_{r},
\end{eqnarray}
with some positive constant $C_{r}$.

Since $\mathcal{E}_{1}$ is continuously embedded into $C^{0}(\mathbb{R},\mathbb{R}^{N})$ with 
 $\| v \|_{\infty} = \sup_{\xi\in \mathbb{R}} |v(\xi)| \le  \| v \|_{1}$, we arrive at the following lemma \cite{Iva16}.
\newtheorem{lemmas}{Lemma}
\begin{lemmas}
$\mathcal{E}_{r}$ is a Hilbert space. 
\end{lemmas}

The next lemma gives \cite{Iva16} an estimate on the absolute value of an element of the space $\mathcal{E}_{r}$.
\begin{lemmas}
If $v\in \mathcal{E}_{r}$ then 
$$| v(\xi) | \le \left(\frac{1+2 p(\xi)}{2 r(\xi)}\right)^{1/2}\| v \|_{r}.$$
\end{lemmas}

Now we may construct a Hilbert manifold modelled on the Hilbert space $\mathcal{E}_{r}$. Consider the set of functions 
\begin{eqnarray}
\nonumber
\mathfrak{M} = \biggl\{q\in AC(\mathbb{R},\mathbb{R}^{N}): q(\xi)\in \mathcal{M} \textrm{ for each } \xi\in\mathbb{R} 
\textrm{ and }& &\\
%\nonumber
 \int\limits_{\mathbb{R}}\biggl( |q'(\xi)|^{2} &+& {\rm d}^{2}\bigl(q(\xi), \chi(\xi)\bigr)  \biggr) 
r(\xi){\rm d}\xi<\infty\biggr\},
\end{eqnarray}
where $d(x, y)$ denotes the Riemannian distance between any $x, y\in \mathcal{M}$
\begin{equation*}
d(x, y) = \inf\limits_{c}\left\{\int\limits_{a}^{b}\vert c'(s) \vert {\rm d}s,\,\,  c: [a,b] \to \mathcal{M} \textrm{ is a piecewise smooth curve}\right\}.
\end{equation*}
and the function $\chi$ is the step-function:
\begin{eqnarray}
\nonumber
\chi(\xi)=
\begin{cases} 
x_{+},\,\, \xi\ge 0,\\
x_{-},\,\, \xi<0.
\end{cases}
\end{eqnarray} 

Then one has the following proposition \cite{Iva16}:
\newtheorem{propositions}{Proposition}
\begin{propositions}
The set $\mathfrak{M}$ is a Hilbert manifold of class $C^{2}$ with tangent space at $q$ given by
\begin{equation}
\nonumber
T_{q}\mathfrak{M} = \biggl\{v\in \mathcal{E}_{r}: v(\xi)\in T_{q(\xi)}\mathcal{M}\,\,  \textrm{ for all } \xi\in\mathbb{R}\biggr\}.
\end{equation}
\end{propositions}

Note that the Lagrangian of the system (1.1) in terms of the new time variable takes the form
$$
L(q, q', \xi) = r(\xi)\left(K(q, q') - \sigma(\xi)V(q)\right).
$$
We introduce new Lagrangian $\hat L(q, q', \xi) = L(q, q', \xi) + r(\xi)\sigma(\xi) V(\chi(\xi))$ and consider the action functional $I$ defined on  $\mathfrak{M}$:
\begin{equation}
I[q] = \int\limits_{\mathbb{R}} \hat L(q, q, \xi){\rm d}\xi.
\end{equation}
It follows from the homogeneity of the quadratic form $K(q, q')$ in velocity $q'$ and the assumption $A_{4}$ (see e.g. \cite{Iva16}) that the integral in (2.2) converges for all $q \in \mathfrak{M}$. Thus, the functional $I$ is well-defined on $\mathfrak{M}$. Moreover, one may prove \cite{Iva16} the following proposition
\begin{propositions}
The functional $I$ is of class $C^{1}(\mathfrak{M})$ with locally Lipschitz derivative. Critical points of $I$ are in one-to-one correspondence with doubly asymptotic trajectories such that $q(\xi) \to x_{\pm}$ and $q'(t)\to 0$ as $\xi\to \pm\infty$.
\end{propositions}
We here note that equations of motion of the system are of the form
\begin{eqnarray}
\frac{d}{d\xi}\left(r(\xi)\frac{\partial T}{\partial q'}\right) - r(\xi)\frac{\partial T}{\partial q} = -r(\xi)\sigma(\xi)\frac{\partial V}{\partial q}.
\end{eqnarray}
and if $q$ is a critical point of $I$ it solves these equations.

The existence of variety of critical points for the functional $I$ follows \cite{Iva16} from the two facts: the first one is fulfillment of the Palais-Smale conditions. One says \cite{PalSm} that a functional $J$ defined on a Hilbert space satisfies the Palais-Smale conditions if any sequence $\{q_{n}\}$ for which $J[q_{n}]$ is bounded and $J'[q_{n}]\to 0$ as $n\to \infty$ possesses a convergent subsequence. It was proved \cite{Iva16} the functional $I$ satisfies the Palais-Smale conditions. Hence the Lusternik-Schnirelmann theory provides a bound on the number of critical points, namely, $\#\{critical points\} \ge cat(\mathfrak{M}) = \infty$ \cite{Serre} and we arrive at the following
\begin{propositions}
The system (1.1)  has infinitely many doubly asymptotic trajectories connecting $x_{-}$ and $x_{+}$.
\end{propositions}

\section{Transversal connecting orbits and non-degenerate critical points of the action functional}

\setcounter{equation}{0}

Let $g$ and $\nabla$ denote a Riemannian metric and the Levi-Civita connection on the manifold $\mathcal{M}$, while $\Gamma=\{\Gamma^{i}_{jk}\}$ stands for the corresponding Christoffel symbols. To emphasize the dependence of the introduced objects on a point $x\in \mathcal{M}$ we will write them as $g(x), \nabla(x), \Gamma(x)$. We also denote by $\langle \cdot, \cdot \rangle$ the scalar product on $T\mathcal{M}$, i.e. for any $x\in \mathcal{M}$ and any $v_{1}, v_{2}\in T_{x}\mathcal{M}$ in local coordinates one has
\begin{equation*}
\langle v_{1}, v_{2}\rangle = g_{ij}(x)v_{1}^{i}v_{2}^{j},
\end{equation*}
where dummy indices summation rule is used. 

Since $\mathcal{M}$ is compact it follows from the Hopf-Rinow theorem that $\mathcal{M}$ is complete in metric $d$ and also geodesically complete, i.e. for any $x, y\in \mathcal{M}$ there exists a geodesic $\varGamma$ connecting $x$ and $y$ with the length $L(\varGamma) = d(x, y)$.
 
Taking this into account the scalar product on $\mathfrak{M}$ denoted by $\langle\langle \cdot, \cdot \rangle\rangle$ reads
\begin{equation}
\langle\langle \varphi_{1}, \varphi_{2}\rangle\rangle = 
\int\limits_{\mathbb{R}}\biggl(
\langle D_{\xi}\varphi_{1}(\xi), D_{\xi}\varphi_{2}(\xi)\rangle +
\langle \varphi_{1}(\xi), \varphi_{2}(\xi)\rangle\biggr)r(\xi) {\rm d}\xi,\,\, 
\forall \varphi_{1}, \varphi_{2}\in T_{q}\mathfrak{M},
\end{equation}
here $D_{\xi}$ stands for the covariant derivative with respect to $\xi$.

Denote by $\mathfrak{X}(\mathcal{M})$ the set of vector fields on the manifold $\mathcal{M}$. One may associate to any tangent vector $v\in T_{x}\mathcal{M}$ a vector field 
$X_{v}\in \mathfrak{X}(\mathcal{M})$ such that $X_{v}(x) = v$. The corresponding covariant differentiation with respect to $X_{v}$ will be denoted by $\nabla_{v}$. Note that the covariant derivative $D_{\xi}$ is related to a vector field $X_{q'(\xi)}$, i.e. $D_{\xi} = \nabla_{q'(\xi)}$.
Besides, the differentiation $\nabla_{v}(x)$ at the point $x$ depends only on the vector $v$, but does not depend on particular representative vector field $X_{v}$.

Then we arrive at the following lemma
\begin{lemmas}
The functional $I$ is of class $C^{2}(\mathfrak{M})$. Moreover,
for any $q\in \mathfrak{M}$ and $\varphi, \psi\in T_{q}\mathfrak{M}$ the first and the second derivatives of $I$ are of the form:
\begin{equation}
I'[q](\varphi) = \int\limits_{\mathbb{R}}
\biggl(\langle D_{\xi}q(\xi), D_{\xi}\varphi(\xi)\rangle - 
\sigma(\xi)\langle {\rm grad}V(q(\xi)), \varphi(\xi)\rangle\biggr)r(\xi){\rm d}\xi,
\end{equation}
\begin{eqnarray}
\nonumber
I''[q](\varphi, \psi) = \int\limits_{\mathbb{R}}
\biggl(\langle D_{\xi}\varphi(\xi), D_{\xi}\psi(\xi)\rangle -
\langle R(D_{\xi}q(\xi), \psi(\xi))\varphi(\xi), D_{\xi}q(\xi)\rangle -\\
\sigma(\xi)\langle H^{V}(q(\xi))\psi(\xi), \varphi(\xi)\rangle\biggr)r(\xi){\rm d}\xi +
I'[q](\nabla_{\psi} X_{\varphi}),
\end{eqnarray}
where $R$ stands for the curvature tensor.
\end{lemmas}
PROOF: - The proof of this lemma is strightforward. Take $q\in \mathfrak{M}$ and $\varphi_{k}\in T_{q}\mathfrak{M}, k=1, 2$. Then one may consider a two-parameter family of proper curves $\alpha$
\begin{equation*}
\alpha: \mathbb{R}\times \prod\limits_{k=1}^{2}\left(-s_{0}^{k}, s_{0}^{k}\right) \to \mathcal{M}
\end{equation*}
defined for some positive $s_{0}^{k}$ by the formula
\begin{equation*}
\alpha(\xi, s_{1}, s_{2}) = 
{\rm Exp}_{q(\xi)}\left(s_{1} \varphi_{1}(\xi)+ s_{2} \varphi_{2}(\xi)\right),
\end{equation*}
where ${\rm Exp}_{x}$ denotes the exponential map at a point $x\in \mathcal{M}$. It satisfies 
\begin{equation*}
X_{\varphi_{k}} = \alpha_{*}\frac{\partial}{\partial s^{k}} \biggr\vert_{s^{1}=s^{2}=0}.
\end{equation*}
Here $\alpha_{*}$ stands for the tangent map to $\alpha$.
\newline
Then one gets
\begin{equation*}
I'[q](\varphi_{k}) = \frac{\partial}{\partial s^{k}} I[\alpha]\biggr\vert_{s^{1}=s^{2}=0},\,\,\,\,\,
I''[q](\varphi_{k}, \varphi_{j}) = \frac{\partial^{2}}{\partial s^{k} \partial s^{j}}   I[\alpha]\biggr\vert_{s^{1}=s^{2}=0}.
\end{equation*}
Define 
\begin{equation*}
X_{k} = \alpha_{*}\frac{\partial}{\partial s^{k}},\,\,\, 
X_{0}=\alpha_{*}\frac{\partial}{\partial \xi}.
\end{equation*}
Then
\begin{align*}
\frac{\partial}{\partial s^{k}} I[\alpha] = \int\limits_{\mathbb{R}}
\Bigl[
\bigl\langle\nabla_{X_{k}}\nabla_{X_{0}}\alpha, \nabla_{X_{0}}\alpha\bigr\rangle - 
\sigma(\xi)\bigl\langle grad V(\alpha), X_{k}\bigr\rangle\Bigr]r(\xi){\rm d}\xi,&\\
\frac{\partial^{2}}{\partial s^{k} \partial s^{j}} I[\alpha] = 
\int\limits_{\mathbb{R}}\Bigl[
\bigl\langle\nabla_{X_{j}}\nabla_{X_{k}}\nabla_{X_{0}}\alpha, \nabla_{X_{0}}\alpha\bigr\rangle +
\bigl\langle\nabla_{X_{k}}\nabla_{X_{0}}\alpha, \nabla_{X_{j}}\nabla_{X_{0}}\alpha\bigr\rangle - &\\
\sigma(\xi)\Bigl(\bigl\langle \nabla_{X_{j}}grad V(\alpha), X_{k}\bigr\rangle+&
\bigl\langle grad V(\alpha), \nabla_{X_{j}}X_{k}\bigr\rangle\Bigr)\Bigr]r(\xi){\rm d}\xi,\\
\end{align*}
where all the integrants are evaluated at a point $(\xi, \vec s)$ with $\vec s = (s_{1}, s_{2})$.

One may note that for all $X, Y, Z\in \mathfrak{X}(\mathcal{M})$ and $1 \le k,j \le 2$ the following equalities hold
\begin{align*}
&X\left(\langle Y, Z\rangle\right) = 
\langle \nabla_{X}Y, Z\rangle + \langle Y,  \nabla_{X}Z\rangle,\quad
\left[X_{k}, X_{0}\right](\xi, \vec s) = 0,
\quad \nabla_{X_{k}}X_{0}(\xi, \vec s) = \nabla_{X_{0}}X_{k}(\xi, \vec s),\\
&\nabla_{X_{j}}\nabla_{X_{0}}X_{k}(\xi, \vec s) = 
\nabla_{X_{0}}\nabla_{X_{j}}X_{k}(\xi, \vec s) + R(X_{j}, X_{0})X_{k}(\xi, \vec s)\\
&\nabla_{X_{k}}V(\alpha(\xi, \vec s)) = 
\langle {\rm grad}V(\alpha(\xi, \vec s)), X_{k}(\xi,\vec s)\rangle,\\
&\langle\nabla_{X_{j}}{\rm grad}V(\alpha(\xi, \vec s)), 
X_{k}(\xi, \vec s)\rangle = 
\langle H^{V}(\alpha(\xi, \vec s))X_{j}(\xi,\vec s), 
X_{k}(\xi,\vec s)\rangle,
\end{align*}
where  $H^{V}$ denotes the Hessian of the function $V$. Taking this into account and setting $s^{1}=s^{2}=0$ proves the desired formulae for the first and the second drivatives of $I$.

We also note that in local coordinates 
\begin{equation*}
\bigl(\nabla_{X_{j}}X_{k}\bigr)^{l} = 
\frac{\partial^{2}\alpha^{l}}{\partial s^{k} \partial s^{j}} +
X_{j}^{s} \Gamma^{l}_{r s}(\alpha) X_{k}^{r},
\end{equation*}
where all summands are evaluated at $(\xi, \vec s)$. $\square$

As a consequence of this lemma we get
\newtheorem{corollaries}{Corollary}
\begin{corollaries}
If $q$ is a doubly asymptotic trajectory connecting $x_{-}$ and $x_{+}$  then for any $\varphi, \psi \in T_{q}\mathfrak{M}$
\begin{equation*}
I''[q](\varphi, \psi) = \int\limits_{\mathbb{R}}
\biggl(\langle D_{\xi}\varphi(\xi), D_{\xi}\psi(\xi)\rangle -
\langle R(D_{\xi}q(\xi), \psi(\xi))\varphi(\xi), D_{\xi}q(\xi)\rangle -\\
\sigma(\xi)\langle H^{V}(q(\xi))\psi(\xi), \varphi(\xi)\rangle\biggr)r(\xi){\rm d}\xi
\end{equation*}
\end{corollaries}

Consider $H^{V}(x_{\pm})$ the Hessian of $V$ at the point $x_{\pm}$. Since $x_{\pm}$ is a maximum point of $\sigma(\xi)V$ and $H^{V}$ is symmetric, one may perform a change of coordinates to diagonalize $\sigma(\xi) H^{V}(x_{\pm}) = - {\rm diag} \left\{{\Lambda_{1}^{\pm}}^{2},...,{\Lambda_{n}^{\pm}}^{2}\right\}$. Further we will refer to such local coordinate system near the point $x_{\pm}$ as to $LC(x_{\pm})$.

Due to definition of the sets $X_{\pm}$ the equilibrium $(x_{+},0)$ (resp. $(x_{-},0)$) possesses an invariant stable (resp. unstable) manifold. One may prove the following proposition \cite{Iva17}
\begin{propositions}
Let $U$ be an open subset of the tangent bundle $T{\mathcal M}$ containing the equilibrium $(x_{+},0)$ and $\Phi_{\xi}$ be the flow of the system $(2.4)$. Then there exists a $n+1$-dimensional differentiable manifold 
$\mathcal{W}^{s}(x_{+})\subset T{\mathcal M}\times \mathbb{R}$ 
 such that for any $\xi_{1} > \xi_{0}$ 
$\Phi_{\xi_{1}-\xi_{0}}(\mathcal{W}^{s}(x_{+},\xi_{0}))\subset \mathcal{W}^{s}(x_{+},\xi_{1})$, where $\mathcal{W}^{s}(x_{+},\xi_{0}) = \{(a,b)\in T\mathcal{M}: (a,b,\xi_{0})\in \mathcal{W}^{s}(x_{+})\}$, and for any 
$(a,b)\in \mathcal{W}^{s}(x_{+},\xi_{0})$
$$\lim\limits_{\xi\to +\infty}\Phi_{\xi}(a,b) = (x_{+},0).$$
 If one replaces the flow $\Phi_{\xi}$ by its inverse and takes the limit $\xi\to -\infty$, the similar statement is valid for the equilibria $(x_{-},0)$ and its unstable manifold $\mathcal{W}^{u}(x_{-})$.

Moreover, if $q(\xi)$ is a solution of $(2.4)$ such that 
$(q(\xi_{0}), q'(\xi_{0}))\in \mathcal{W}^{s}(x_{+},\xi_{0})$ 
(resp. $\mathcal{W}^{u}(x_{-},\xi_{0})$) then for any $\lambda^{\pm}$ satisfying
$$0<\lambda^{\pm} < \Lambda^{\pm}_{\min},\quad 
\Lambda^{\pm}_{\min} = \min\limits_{k=1,\ldots,n}\left\{\Lambda_{k}^{\pm}\right\}$$
the following estimate holds
$$d(q(\xi), x_{\pm}) =  
O\left(\frac{1}{r^{1/4}(\xi)} 
{\rm e}^{\mp\lambda^{\pm} \int\limits_{0}^{\xi}r(s){\rm d}s}\right),\quad \xi\to \pm\infty.$$
\end{propositions}

Let $q$ be a doubly asymptotic trajectory which existence is guaranteed by Proposition 3. Then for any $\xi\in \mathbb{R}$ $(q(\xi), q'(\xi)) \in \mathcal{W}^{s}(x_{+},\xi)\cap \mathcal{W}^{u}(x_{-},\xi)$.
\newline
We say that the trajectory $q$ is {\it transversal} if for any $\xi\in \mathbb{R}$ $\mathcal{W}^{s}(x_{+},\xi)$ intersects $\mathcal{W}^{u}(x_{-},\xi)$ transversally at the point $(q(\xi), q'(\xi))$.

One may characterize such transversal connecting orbits in a different way:

\begin{propositions}
A doubly asymptotic trajectory $q$ is transversal if and only if $q$ is non-degenerate critical point of the action functional $I$.
\end{propositions}
PROOF: - First we note that if $q$ is a critical point of the functional $I$ it is a solution of the equations of motion. Hence,
if $\mathcal{W}^{s}(x_{+},\xi)$ and $\mathcal{W}^{u}(x_{-},\xi)$ intersect transversally at a point $(q(\xi), q'(\xi))$ for some $\xi\in \mathbb{R}$ they do so for any $\xi\in \mathbb{R}$.

Represent the second derivative of $I$ at $q$ as 
\begin{equation*}
I''[q](\varphi, \psi) = \langle \mathfrak{A} \varphi, \psi\rangle_{\mathbb{L}_{2}},
\end{equation*}
where $\varphi, \psi\in T_{q}\mathfrak{M}$, $\mathfrak{A}$ is a differential operator 
\begin{equation}
\mathfrak{A} = -D_{\xi}r(\xi)D_{\xi} - r(\xi)\left(R(D_{\xi}, q'(\xi)) q'(\xi) + \sigma(\xi)H^{V}(q(\xi))\right) 
\end{equation}
and $\langle \varphi, \psi\rangle_{\mathbb{L}_{2}} = \int\limits_{\mathbb{R}}\varphi(\xi)\psi(\xi){\rm d}t$.
The nullspace of $\mathfrak{A}$ consists of those $\varphi\in T_{q}\mathfrak{M}$ which solve
\begin{equation}
D_{\xi}r(\xi)D_{\xi}\varphi + r(\xi)\left(R(D_{\xi}\varphi, q'(\xi)) q'(\xi) + 
\sigma(\xi)H^{V}(q(\xi))\varphi\right) = 0.
\end{equation}
Note that equation (3.5) is the variational equation of (2.4) along the trajectory $q$. It was proved in \cite{Iva17} that this equation possesses exponential dichotomy on $\mathbb{R}_{\pm}$. Moreover, if $(\varphi_{0}, {\varphi'}_{0})\notin T_{(q(\xi_{0}), q'(\xi_{0}))}\mathcal{W}^{s}(x_{+}, \xi_{0})$ (resp. $(\varphi_{0}, {\varphi'}_{0})\notin T_{(q(\xi_{0}), q'(\xi_{0}))}\mathcal{W}^{u}(x_{-}, \xi_{0})$), the solution $\varphi(\xi) = \Phi_{0}(\xi - \xi_{0}) (\varphi_{0}, {\varphi'}_{0}, \xi_{0})$ satisfies
\begin{equation}
\vert \varphi(\xi) \vert > C(\varphi_{0}, {\varphi'}_{0}) \frac{1}{r^{1/4}(\xi)} 
{\rm e}^{\pm\lambda^{\pm} \int\limits_{0}^{\xi}r(s){\rm d}s},\quad \xi\to \pm\infty,
\end{equation}
where $\lambda^{\pm}$ is an arbitrary constant such that $0<\lambda^{\pm}<\Lambda_{\min}^{\pm}$ and $\Lambda_{\min}^{\pm}$ is defined in Proposition 4.
Hence, for any point $(q(\xi), q'(\xi))$ one has 
\begin{equation*}
T_{(q(\xi), q'(\xi))}T\mathcal{M} = T_{(q(\xi), q'(\xi))}\mathcal{W}^{s}(x_{+}, \xi)\oplus \mathcal{N}^{u}(x_{+},\xi) = 
T_{(q(\xi), q'(\xi))}\mathcal{W}^{u}(x_{-}, \xi)\oplus \mathcal{N}^{s}(x_{-},\xi),
\end{equation*}
where $\mathcal{N}^{u}(x_{+},\xi), \mathcal{N}^{s}(x_{-},\xi)$ are orthogonal complements to 
$T_{(q(\xi), q'(\xi))}\mathcal{W}^{s}(x_{+}, \xi)$ and $T_{(q(\xi), q'(\xi))}\mathcal{W}^{u}(x_{-}, \xi)$, respectively. Besides, if 
$(\varphi_{0}, {\varphi'}_{0})\in \mathcal{N}^{u}(x_{+}, \xi_{0})$ 
(resp. $\mathcal{N}^{s}(x_{-}, \xi_{0}))$, the estimate (3.6) holds. It means that there exist exactly $n$ linear independent solutions $\varphi^{+}_{k}, k=1,\ldots, n$ of (3.5) which decay at $+\infty$ and $n$ linear independent solutions $\varphi^{-}_{k}, k=1,\ldots, n$ of (3.5) which decay at $-\infty$. 

If $q$ is transversal then $\mathcal{W}^{s}(x_{+}, \xi)$ and $\mathcal{W}^{u}(x_{-}, \xi)$ are transverse at $(q(\xi), q'(\xi))$ and, hence, the solutions $\varphi^{\pm}_{k}, k=1,\ldots, n$ are linearly independent. It implies that the nullspace of $\mathfrak{A}$ is trivial and due to Fredholmness of $\mathfrak{A}$ its invertibility. Thus $q$ is non-degenerate critical point.
Since all implications are valid in both directions the opposite statement follows. $\square$

\section{Newton-Kantorovich theorem}

\setcounter{equation}{0}

In this section we collect some results on the Newton's method.

The manifold $\mathfrak{M}$ is an infinite dimensional Hilbert manifold which is topologically equivalent to the space of paths 
$\Omega(\mathcal{M}, x_{-}, x_{+})$ connecting $x_{-}$ to $x_{+}$. As in finite dimensional case one may construct a unique, symmetric connection (the Levi-Civita connection), compatible with the Riemannian metric (3.1) \cite{Lang99}. This leads to definition of parallel transport, geodesics, exponential map, the curvature tensor etc. To distinct the mentioned objects from the corresponding objects of the manifold $\mathcal{M}$ the former will be supplemented by the sign '$\hat{\phantom{x}}$' (e.g. the Levi-Civita connection on $\mathfrak{M}$ will be denoted by $\hat \nabla$).

Let $X, Y, Z$ be $C^{1}$-vector fields on $\mathfrak{M}$ then the Levi-Civita connection $\hat \nabla$ at a point 
$q\in \mathfrak{M}$ is given \cite{MaeRosTor} by
\begin{align*}
%\nonumber
\langle\langle \hat\nabla_{X}Y, Z\rangle\rangle = 
&\int\limits_{\mathbb{R}}\biggl( \langle \nabla_{q'(\xi)} \nabla_{X(\xi)}Y(\xi), \nabla_{q'(\xi)} Z(\xi)\rangle + 
\langle \nabla_{X(\xi)}Y(\xi), Z(\xi)\rangle -\\
%\nonumber
&-\frac{1}{2}\biggl[\langle \nabla_{q'(\xi)}(R(X(\xi), q'(\xi))Y(\xi)), Z(\xi)\rangle + 
\langle R(X(\xi), q'(\xi))\nabla_{q'(\xi)} Y(\xi), Z(\xi)\rangle +\\
%\nonumber
&\langle \nabla_{q'(\xi)}(R(Y(\xi), q'(\xi))X(\xi)), Z(\xi)\rangle + 
\langle R(Y(\xi), q'(\xi))\nabla_{q'(\xi)} X(\xi), Z(\xi)\rangle +\\
&\langle R(\nabla_{q'(\xi)} X(\xi), Y(\xi)) q(\xi), Z(\xi)\rangle -
\langle R(\nabla_{q'(\xi)} X(\xi), Y(\xi)) q'(t), Z(\xi)\rangle\biggr]\biggr)r(\xi){\rm d}\xi. 
\end{align*}

For any $X\in \mathfrak{X}(\mathfrak{M})$ the covariant derivative at a point $q\in \mathfrak{M}$ defines a linear map
\begin{eqnarray}
\nonumber
\mathcal{D} X(q): T_{q}\mathfrak{M}\to  T_{q}\mathfrak{M},\\
\nonumber
\varphi \mapsto \hat\nabla_{(\varphi)} X(q).
\end{eqnarray}

Let $a, b\in \mathbb{R}$ and $\alpha: (a, b) \to \mathfrak{M}$ be a piecewise smooth curve. We fix $s_{0}\in (a, b)$ and denote the parallel transport along $\alpha$ by 
$\hat P_{\alpha, s_{0}, s}$. Then
\begin{eqnarray}
\nonumber
\hat P_{\alpha, s_{0}, s}: T_{\alpha(s_{0})}\mathfrak{M} \to T_{\alpha(s)}\mathfrak{M},\\
\nonumber
\varphi \mapsto \hat P_{\alpha, s_{0}, s}(\varphi) = \gamma(s; \varphi),
\end{eqnarray}
where $\gamma: (a, b)\to T\mathfrak{M}$ is the unique curve which is $\alpha$-parallel and $\gamma(s_{0})=\varphi$. Note that $\hat P_{\alpha, s_{0}, s}$ is linear and satisfies for all $s_{0}, s_{1}, s_{2}\in (a, b)$
\begin{eqnarray*}
\nonumber
\hat P_{\alpha, s_{1}, s_{2}}\circ \hat P_{\alpha, s_{0}, s_{1}} = \hat P_{\alpha, s_{0}, s_{2}},\\
\hat P_{\alpha, s_{1}, s_{0}} = \hat P_{\alpha, s_{0}, s_{1}}^{-1}.
\end{eqnarray*}

If $X$ is a vector field on $\mathfrak{M}$, then one has
\begin{equation}
\mathcal{D}X(\alpha(s)) \alpha'(s) = \hat\nabla_{\alpha'(s)} X(\alpha(s)) = 
\lim\limits_{h\to 0}\frac{1}{h}\left( \hat P_{\alpha, s+h, s} X(\alpha(s+h)) - X(\alpha(s))\right).
\end{equation}
Indeed, to prove the second equality rewrite it in a local chart $\mathfrak{U}$ around $\alpha(s)$:
\begin{equation}
\hat\nabla_{\alpha'_{\mathfrak{U}}(s)} X_{\mathfrak{U}}(\alpha_{\mathfrak{U}}(s)) = 
\lim\limits_{h\to 0}\frac{1}{h}\left( \hat P^{\mathfrak{U}}_{\alpha, s+h, s}
X_{\mathfrak{U}}(\alpha_{\mathfrak{U}}(s+h)) - 
X_{\mathfrak{U}}(\alpha_{\mathfrak{U}}(s))\right),
\end{equation}
where the index $\mathfrak{U}$ indicates the corresponding representatives in the chart. We also note (see e.g. \cite{Lang99}) that for any two vector fields $X, Y\in \mathfrak{X}(\mathfrak{M})$ 
\begin{equation}
(\hat\nabla_{X}Y)_{\mathfrak{U}}(q) = Y'_{\mathfrak{U}}(q) \cdot X_{\mathfrak{U}}(q) - 
B_{\mathfrak{U}}(q; X_{\mathfrak{U}}(q), Y_{\mathfrak{U}}(q)),
\end{equation}
where $Y'_{\mathfrak{U}}(q)$ is the tangent map to the section 
$Y_{\mathfrak{U}}: \mathfrak{U}\to \mathfrak{U}\times \mathcal{E}_{r}$ and 
$B_{\mathfrak{U}}(q; u, v)$ is the local representative of the symmetric bilinear map associated with the unique metric spray.

Then, using linearity of the parallel transport and formula for derivative of a composition, one obtains
\begin{align}
\nonumber
&\lim\limits_{h\to 0}\frac{1}{h}\left( \hat P^{\mathfrak{U}}_{\alpha, s+h, s} 
X_{\mathfrak{U}}(\alpha_{\mathfrak{U}}(s+h)) - 
X_{\mathfrak{U}}(\alpha_{\mathfrak{U}}(s))\right)=\\
\nonumber
&\lim\limits_{h\to 0}\frac{1}{h}\hat P^{\mathfrak{U}}_{\alpha, s+h, s} 
\left(X_{\mathfrak{U}}(\alpha_{\mathfrak{U}}(s+h)) - 
\hat P^{\mathfrak{U}}_{\alpha, s, s+h} X_{\mathfrak{U}}(\alpha_{\mathfrak{U}}(s))\right)=\\
\nonumber
&\lim\limits_{h\to 0}\frac{1}{h}\hat P^{\mathfrak{U}}_{\alpha, s+h, s} 
\left(X_{\mathfrak{U}}(\alpha_{\mathfrak{U}}(s+h)) - 
X_{\mathfrak{U}}(\alpha_{\mathfrak{U}}(s)) +
X_{\mathfrak{U}}(\alpha_{\mathfrak{U}}(s)) -
\hat P^{\mathfrak{U}}_{\alpha, s, s+h} X_{\mathfrak{U}}(\alpha_{\mathfrak{U}}(s))\right)=\\
\nonumber
&\lim\limits_{h\to 0}\frac{1}{h}\left[
\hat P^{\mathfrak{U}}_{\alpha, s+h, s} \left(X_{\mathfrak{U}}(\alpha_{\mathfrak{U}}(s+h)) - X_{\mathfrak{U}}(\alpha_{\mathfrak{U}}(s)) +
\beta_{\mathfrak{U}}(s; X_{\mathfrak{U}}(\alpha_{\mathfrak{U}}(s))) -
\beta_{\mathfrak{U}}(s+h; X_{\mathfrak{U}}(\alpha_{\mathfrak{U}}(s)))\right)\right]=\\
&X'_{\mathfrak{U}}(\alpha_{\mathfrak{U}}(s))\cdot \alpha'_{\mathfrak{U}}(s) - 
\beta'_{\mathfrak{U}}(s; X_{\mathfrak{U}}(\alpha_{\mathfrak{U}}(s))),
\end{align}
where $\beta_{\mathfrak{U}}(s; X_{\mathfrak{U}}(\alpha_{\mathfrak{U}}(s_{0}))$ stands for $\alpha$-parallel curve in $\mathfrak{U}\times \mathcal{E}_{r}$ such that
$\beta'_{\mathfrak{U}}(s_{0}; X_{\mathfrak{U}}(\alpha_{\mathfrak{U}}(s_{0}))) = 
X_{\mathfrak{U}}(\alpha_{\mathfrak{U}}(s_{0})).$

Since $\beta_{\mathfrak{U}}$ is $\alpha$-parallel and due to (4.3) we have 
\begin{equation}
\beta'_{\mathfrak{U}}(s; X_{\mathfrak{U}}(\alpha_{\mathfrak{U}}(s))) = 
B_{\mathfrak{U}}(\alpha_{\mathfrak{U}}(s); \alpha'_{\mathfrak{U}}(s), X_{\mathfrak{U}}(\alpha_{\mathfrak{U}}(s))).
\end{equation}
Substituting this into (4.4) and using (4.3) together with symmetry of $B_{\mathfrak{U}}$ we prove (4.2).

Let $\mathfrak{U}$ be an open subset of $\mathfrak{M}$, $X\in \mathfrak{X}(\mathfrak{U})$. We will say that the covariant derivative 
$\mathcal{D}X$ is Lipschitz with constant $C_{L}>0$ 
($\mathcal{D}X\in {\rm Lip}_{C_{L}}(\mathfrak{U})$) if for any geodesic $\alpha$ and 
$a, b\in \mathbb{R}$ such that $\alpha([a, b])\subset \mathfrak{U}$
\begin{equation}
\Vert \hat P_{\alpha, b, a} \mathcal{D}X(\alpha(b)) \hat P_{\alpha, a, b} - \mathcal{D}X(\alpha(a))\Vert_{op}\ge
C_{L}\int\limits_{a}^{b} \Vert \alpha'(s)\Vert {\rm d}s,
\end{equation}
where $\Vert\cdot \Vert_{op}$ stands for the operator norm.

Then one has the following lemma \cite{FerSva}
\begin{lemmas}
Let $\mathfrak{U}$ be an open subset of $\mathfrak{M}$ and $X$ is a vector field being continuous on $\overline{\mathfrak{U}}$ and $C^{1}$ on $\mathfrak{U}$ with 
$\mathcal{D}X\in {\rm Lip}_{C_{L}}(\mathfrak{U})$. For given $q\in \mathfrak{U}$ and $\varphi\in T_{q}\mathfrak{M}$ define a geodesic
$$\alpha(s) = \hat{\rm Exp}_{q}(s \varphi).$$
If $\alpha([0,s))\subseteq \mathfrak{U}$ then
$$\hat P_{\alpha, s, 0}X(\alpha(s)) = X(q) + s \mathcal{D}X(q) \varphi + Rem(s),\,\,\, 
\Vert Rem(s)\Vert \le \frac{C_{L}}{2} s^{2} \Vert \varphi \Vert^{2}.$$
\end{lemmas}

We introduce the following notations for an open and closed ball in $\mathfrak{M}$:
\begin{align*}
%\nonumber
&B(q_{0}, r) = \{ q\in \mathfrak{M}: \hat d(q, q_{0}) < r\},\\
%\nonumber
&B[q_{0}, r] = \{ q\in \mathfrak{M}: \hat d(q, q_{0}) \le r\}.
\end{align*}

Assume the conditions of Lemma 4 are fulfilled. Take any point $q_{0}\in U$ and define the so-called Newton's sequence
\begin{equation}
q_{k+1} = \hat {\rm Exp}_{q_{k}}\left(- \mathcal{D}X(q_{k})^{-1}X(q_{k})\right).
\end{equation}
The following theorem gives conditions which guarantee the well definedness of the Newton's sequence and its convergence to a singular point of the vector field $X$.
\newtheorem{theorems}{Theorem}
\begin{theorems}[Kantorovich's theorem in Riemannian manifold, \cite{FerSva}]
Let $\mathfrak{U}$ be an open subset of $\mathfrak{M}$ and $X$ is a vector field being continuous on $\overline{\mathfrak{U}}$ and $C^{1}$ on $\mathfrak{U}$ with $\mathcal{D}X\in {\rm Lip}_{C_{L}}(\mathfrak{U})$. Let $q_{0}\in \mathfrak{U}$ such that $\mathcal{D}X(q_{0})$ is nonsingular and for some positive $a, b \in \mathbb{R}$
 $$\Vert \mathcal{D}X(q_{0})^{-1}\Vert_{op}\le a,\,\,\,  
\Vert \mathcal{D}X(q_{0})^{-1}X(q_{0})\Vert \le b,\,\,\,
l = a b C_{L}\le 1/2$$
and
$$B(q_{0}, r_{*})\subseteq \mathfrak{U},$$
where $r_{*} = \frac{1}{aC_{L}}\left(1 - \sqrt{1 - 2l}\right)$. Then the Newton's sequence $\{q_{k}\}$ generated by the starting point $q_{0}$ is well defined and contained in $B(q_{0}, r_{*})$ and converges to a point $q_{*}$ which is the unique singularity of $X$ in $B[q_{0}, r_{*}]$. Moreover, the rate of convergence is quadratic:
$$\hat d(q_{k}, q_{*})\le (2 l)^{2^{k}}\frac{b}{l},\,\,\, k = 1,2,\ldots .$$
\end{theorems}
{\bf Remark}
We note here that Theorem 1 was proved in \cite{FerSva} for connected geodesically complete Riemannian manifolds which are of finite dimension. However, as the authors of \cite{FerSva} emphasized the finite dimensionality plays no role in the original Kantorovich theorem for Banach spaces \cite{Kan}. One may point out that finite-dimensionality appears in the proof of Theorem 1 of \cite{FerSva} only in two places. The first one is the formula (4.1). It was proved in \cite{FerSva} only under finite dimesion assumption. However, as was shown earlier it is also valid for infinite dimension. The second place is the assumption that every two point of the manifold can be joined by a geodesic. This is much more restrictive. Due to Hopf-Rinow theorem the metric and geodesic completnesses are equivalent in finite-dimensional case. Moreover, in that case it is equivalent to the fact that any two points can be connected by a geodesic. In the infinite dimensional case this is not true. We may only guarantee the metric completeness implies the geodesic one. Although due to result of Ekeland \cite{Ekel} the Hopf-Rinow theorem is generically satisfied, i.e. for any point $q$ of the complete Hilbert manifold the set of points $p$ which cannot be joined by a minimal geodesic to $q$ is of the first category and in particular is nowhere dense. Nevertheless, counterexamples of Grossman \cite{Gros}, McAlpin and Atkin \cite{Atk} show that one may construct a Hilbert manifold which is metrically and geodiesicallly complete, but with points which cannot be joined not only by a minimal geodesic but by a geodesic at all. Thus, the exponential map need not be surjective in the infinite-dimensional case. It is to be noted that existence of a minimal geodesic used in the proof of Theorem 1 in \cite{FerSva} only locally, namely, in a ball of small radius centered at a point $q$ and one might hope to apply local invertability of the exponential map which is valid for metrically complete manifolds. However, the size of a neighborhood where the exponential map is diffeomorphism often is not known. A complete Hilbert manifold is called Hopf-Rinow manifold if any two points of this manifold can be joined by a minimal geodesic. There are different examples of Hopf-Rinow manifolds. The most interesting for a purpose of this paper example is due to Eliason \cite{Eli} which shows that the Sobolev manifolds, i.e. the spaces of the Sobolev sections of a vector bundle on a compact manifold, are Hopf-Rinow. One may conclude that if a manifold in conditions of Theorem 1 is Hopf-Rinow then all statements of this theorem are valid without any changes in the proof.

Using Riemannian structure (3.1) on the manifold $\mathfrak{M}$ we define a vector field $\hat X_{I}$ such that for any 
$q\in \mathfrak{M}$ and any $\varphi\in T_{q}\mathfrak{M}$
\begin{equation}
I'[q](\varphi) = \langle\langle \hat X_{I}(q), \varphi \rangle\rangle.
\end{equation}
Then the doubly asymptotic trajectories connecting $x_{-}$ and $x_{+}$ correspond to singular points of $\hat X_{I}$ and one may try to apply the Newton's method for finding such points. 

\section{Algorithm for construction of doubly asymptotic trajectories}

\setcounter{equation}{0}

In this section we describe a procedure which allows constructing transverse connecting orbits for the system (1.1).

 Let $\{\xi_{k}\}_{k=0}^{\infty}$ be an increasing sequence of positive real numbers and $\Omega_{k} = [-\xi_{k}, \xi_{k}]$. For each interval $\Omega_{k}$ we define the following objects:
\begin{eqnarray}
\nonumber
\mathcal{E}_{r, k} = \biggl\{v\in AC(\Omega_{k},\mathbb{R}^{N}): 
v(\pm\xi_{k}) = 0,\,\,\,
\|v\|_{r,k}^{2} = \int\limits_{\Omega_{k}} \left( |v'(\xi)|^{2} +  |v(\xi)|^{2}\right)
r(\xi){\rm d}\xi < \infty\biggr\},\\
\nonumber
\mathfrak{M}_{k} = \biggl\{q\in AC(\Omega_{k},\mathbb{R}^{N}): q(\xi)\in \mathcal{M} \textrm{ for each }\xi \in\Omega_{k},\,\,\,
q(\pm\xi_{k}) = x_{\pm}\,\,\,
\textrm{ and }\\
\nonumber
 \int\limits_{\Omega_{k}}\biggl( |q'(\xi)|^{2} + d^{2}(q(\xi), \chi(\xi)) \biggr) 
r(\xi){\rm d}\xi<\infty\biggr\},
\end{eqnarray}

Then we arrive at the following proposition
\begin{propositions}
For any integer $k\ge 0$ the set $\mathcal{E}_{r,k}$ is a Hilbert space. If $v\in \mathcal{E}_{r,k}$ then 
$$| v(\xi) | \le \left(\frac{1+2 p(\xi)}{2 r(\xi)}\right)^{1/2}\| v \|_{r,k}.$$
The set $\mathfrak{M}_{k}$ is a Hilbert manifold of class $C^{2}$ with tangent space at $q$ given by
\begin{equation}
\nonumber
T_{q}\mathfrak{M}_{k} = \biggl\{v\in \mathcal{E}_{r,k}: v(\xi)\in T_{q(\xi)}\mathcal{M}\,\,  \textrm{ for all } \xi\in\Omega_{k}\biggr\}.
\end{equation}
\end{propositions}
PROOF:- The proof of this proposition repeats the proof of corresponding statements of section 2.

Note that if $q\in \mathfrak{M}_{k}$ and $v\in T_{q}\mathfrak{M}_{k}$ one may define
\begin{equation*}
\tilde q(\xi) = 
\begin{cases}
q(\xi),\,\,\, \xi\in \Omega_{k},\\
\chi(\xi),\,\,\, \xi\in \Omega_{k+1}\setminus \Omega_{k}
\end{cases},\quad
\tilde v(\xi) = 
\begin{cases}
v(\xi),\,\,\, \xi\in \Omega_{k},\\
0,\,\,\, \xi\in \Omega_{k+1}\setminus \Omega_{k}
\end{cases}.
\end{equation*}
Then $\tilde q\in \mathfrak{M}_{k+1}$ and $\tilde v\in T_{\tilde q}\mathfrak{M}_{k+1}$. Using this identification we may write $\mathfrak{M}_{k}\subset \mathfrak{M}_{k+1}$ and $T_{q}\mathfrak{M}_{k}\subset T_{q}\mathfrak{M}_{k+1}$. However, note that $T_{q}\mathfrak{M}_{k}\neq T_{\tilde q}\mathfrak{M}_{k+1}$. For sake of simplicity in the rest of the paper we will not distinguish between $q$ and its continuation $\tilde q$.

Denote by $Q_{k}(\mathcal{M}, x_{-}, x_{+})$ the set of solutions $q(\xi)$ of (2.4) satisfying $q(\pm \xi_{k}) = x_{\pm}$.

Now we consider the action functional  $I_{k}$ defined on $\mathfrak{M}_{k}$ by the formula 
\begin{equation}
I_{k}[q] = \int\limits_{\Omega_{k}} \hat L(q, q', \xi){\rm d}\xi.
\end{equation}
Then as for the case of the whole real line one may prove the following
\begin{propositions}
The functional $I_{k}$ is of class $C^{2}(\mathfrak{M}_{k})$. The set of critical points of $I_{k}$ coincides with $Q_{k}(\mathcal{M}, x_{-}, x_{+})$. Moreover, the first and the second derivatives of $I_{k}$ are defined by the formulae (3.2), (3.3), where integration over $\mathbb{R}$ is replaced by integration over $\Omega_{k}$.
\end{propositions}

Assume we have constructed a solution $q_{k}\in Q_{k}(\mathcal{M}, x_{-}, x_{+})$ which is a non-singular critical point of $I_{k}$. Using $q_{k}$ as initial approximation one may try to apply the Newton-Kantorovich theorem to construct a solution  
$q_{k+1}\in Q_{k+1}(\mathcal{M}, x_{-}, x_{+})$. As it follows from the Theorem 1 to perform this "k+1"-step one needs to control the smallness of $I'_{k+1}[q_{k}]$ and non-degeneracy of $I''_{k+1}[q_{k}]$.

To estimate the derivatives of the action functional we consider the following auxiliary boundary value problem defined on an interval $[a, b]$ with some $0 < a < b$:
\begin{eqnarray}
\nonumber
\frac{1}{r(\xi)}\frac{{\rm d}}{{\rm d}\xi}r(\xi) \frac{{\rm d}}{{\rm d}\xi} v(\xi) = 
\lambda^{2} v(\xi), \\
\nonumber
v(a) = 1,\,\, v(b) = 0
\end{eqnarray}
or equivalently
\begin{eqnarray}
\nonumber
v'' + 2p(\xi)v' - \lambda^{2} v(\xi) = 0, \\
v(a) = 1,\,\, v(b) = 0.
\end{eqnarray}

The function $p$ is non-negative and decreasing. Hence,
$0\le p(b) = p_{-} \le p(\xi)\le p_{+} = p(a)$. Denote by $v_{\pm}(\xi)$ the solutions of (5.2) where $p(\xi)$ is replaced by $p_{\pm}$, respectively. Then we arrive at the following
\begin{lemmas}
The unique solution $v$ of the boudary value problem (5.2) is strictly decreasing on the interval $[a, b]$ and satisfies for all $\xi\in [a, b]$ to
\begin{eqnarray}
v_{+}(\xi)\le v(\xi)\le v_{-}(\xi)
\end{eqnarray}
\end{lemmas}
PROOF:- Introduce
\begin{equation*}
%\nonumber
\theta_{\pm}(\zeta)={\rm e}^{-p_{\pm} \zeta}\sinh\left(\sqrt{\lambda^{2}+p_{\pm}^{2}}\zeta\right),\quad
\zeta =b - \xi.
\end{equation*}
Then $v_{\pm}(\xi) = \theta_{\pm}(b-\xi)/\theta_{\pm}(b-a)$ and
\begin{eqnarray}
\nonumber
v(\xi) = v_{\pm}(\xi) + \frac{\theta_{\pm}(b-\xi)}{\theta'_{\pm}(0)}
\left[
F_{\pm}(b-\xi) - F_{\pm}(b-a)
%\theta_{\pm}^{-1}(\xi)\int\limits_{0}^{\xi}\theta_{\pm}(\xi-s)F_{\pm}(s){\rm d}s -
%\theta_{\pm}^{-1}(\xi_{0})\int\limits_{0}^{\xi_{0}}\theta_{\pm}(\xi_{0}-s)F_{\pm}(s){\rm d}s
\right],
\end{eqnarray}
where 
\begin{eqnarray}
\nonumber
F_{\pm}(\zeta) = 
2\theta_{\pm}^{-1}(\zeta)\int\limits_{0}^{\zeta}\theta_{\pm}(\zeta-s)
(p(b-s) - p_{\pm})v'(b-s){\rm d}s.
\end{eqnarray}
Hence
\begin{eqnarray}
\nonumber
F'_{\pm}(\zeta) = 
2\theta_{\pm}^{-1}(\zeta)\int\limits_{0}^{\zeta}\left[
\frac{\theta'_{\pm}(\zeta-s)}{\theta_{\pm}(\zeta-s)} -
\frac{\theta'_{\pm}(\zeta)}{\theta_{\pm}(\zeta)}\right]
\theta_{\pm}(\zeta-s)
(p(b-s) - p_{\pm})v'(b-s){\rm d}s.
\end{eqnarray}
Note that $\theta_{\pm}(\zeta)\ge 0$ and $\theta'_{\pm}(\zeta)/\theta_{\pm}(\zeta)$ decreases on $(0, b-a]$. Moreover, by the Shturm's comparison lemma it follows that
$\exp\left(\int\limits_{a}^{b}p(s){\rm d}s\right)v(\xi)v'(\xi)$ increases on $[a, b]$. Since $v(b)v'(b) = 0$ then $v(\xi)v'(\xi)<0$ for all $\xi\in [a, b]$. Consequently, $v(\xi)$, $-v'(\xi)$ are positive on $[a, b)$. Taking all this into account we conclude that $F'_{+}(\zeta)\ge 0$ and $F'_{-}(\zeta)\le 0$ for all $\zeta\in [0, b-a]$. This finishes the proof. $\square$

Now we consider a modified boundary value problem:
\begin{eqnarray}
\nonumber
\frac{1}{r(\xi)}\frac{{\rm d}}{{\rm d}\xi}r(\xi) \frac{{\rm d}}{{\rm d}\xi} v(\xi) = 
F(v(\xi), \xi), \\
v(a) = 1,\,\, v(b) = 0,
\end{eqnarray}
where $F(v,\xi)$ is a continuous function such that 
\begin{align} 
%\nonumber
G(v,\xi) = \int\limits_{0}^{v}F(p, \xi){\rm d}p > \frac{\lambda^{2}}{2}v^{2} \quad \forall\,\,  \xi \in [a, b].
\end{align}
\begin{lemmas}
The unique solution $v$ of the boudary value problem (5.4) is strictly decreasing on the interval $[a, b]$. Moreover,  for any two functions $F_{j}, j=1,2$ satisfying (5.5) and 
$F_{2}(v, \xi) > F_{1}(v, \xi)$ for all $\xi\in [a, b]$, the corresponding solutions $v_{j}, j=1,2$ of the problem (5.4) satisfy
$$
v_{2}(\xi) < v_{1}(\xi)\quad \forall\,\, \xi\in (a, b).
$$
\end{lemmas}
PROOF:- First we observe that if there exists $\xi_{*}\in (a, b)$ such that $v'(\xi_{*}) = 0$ and $v(\xi_{*})>0$ then due to $v(b)=0$ there exists $\xi_{**}\in (\xi_{*}, b)$ such that 
$(r(\xi) v'(\xi))'\vert_{\xi_{**}} < 0$ what contradicts to (5.5) (this follows from equality 
$ r(\xi)\partial_{\xi} = \partial_{t}$ and $\partial_{t}^{2} v = F(v, \xi(t)) > 0$). In the case $v(\xi_{*})\le 0$ there exists $\xi_{**}\in (a, \xi_{*}]$ such that $v(\xi_{**})=0$. However, for each interval $[c, d]\subset [a, b]$ the solution $v$ minimizes the functional 
\begin{equation*}
I_{F}[w] = \int\limits_{c}^{d}\left(\frac{1}{2}\vert w'(\xi)\vert^{2} + G(w(\xi),\xi)\right)r(\xi){\rm d}\xi
\end{equation*}
defined on the set of absolutely continuous functions satisfying boundary conditions $w(c) = v(c), w(d) = v(d)$. But on the interval $[\xi_{**}, b]$ the functional $I_{F}$ attains its minimum at $w=0$. Hence, $v(\xi)\equiv 0$ on $[\xi_{**}, b]$ what contradicts to (5.5).

To prove the second statement of the lemma we consider $w = v_{2} - v_{1}$. It satisfies
\begin{align*}
%\nonumber
&\frac{1}{r(\xi)}\frac{{\rm d}}{{\rm d}\xi}r(\xi) \frac{{\rm d}}{{\rm d}\xi} w(\xi) = 
F_{2}(v_{2}(\xi), \xi) - F_{1}(v_{1}(\xi), \xi), \\
&w(a) = 0,\,\, w(b) = 0.
\end{align*}
If one assumes the existence of $\xi_{*}\in (a, b)$ such that $w'(\xi_{*}) = 0$ and $w(\xi_{*})\ge 0$ then there should exist $\xi_{**}\in [\xi_{*}, b)$ such that 
$(r(\xi) w'(\xi))'\vert_{\xi_{**}} < 0$ what contradicts to $F_{2}(v_{2}(\xi_{**}), \xi_{**}) > F_{1}(v_{1}(\xi_{**}), \xi_{**})$. Hence, $w(\xi)<0\,\, \forall\,\, \xi\in (a, b)$. $\square$

As a collorary to these lemmae we get
\newtheorem{colloraries}{Collorary}
\begin{colloraries}
Let $v$ be the solution of the boudary value problem (5.4) with $F(v, \xi)$ satisfying (5.5) then
$$
\vert v'(b)\vert < \frac{\lambda}{\sinh\left(\lambda (b-a)\right)}.
$$
\end{colloraries}

For any fixed $\lambda$ satisfying $0 < \lambda < \min\{\Lambda_{k}^{\pm}, k = 1,\ldots, n\}$ we consider $R_{\lambda}>0$ such that there exist local charts $U^{\pm}$ around $x_{\pm}$, 
$B_{R_{\lambda}}(x_{\pm})\subset U^{\pm}$ and for any $x\in B_{R_{\lambda}}(x_{\pm})$ and $v\in T_{x}\mathcal{M}$ the following estimates take place:
\begin{align}
\nonumber
V(x_{\pm}) - V(x) &\ge \frac{1}{2}\lambda^{2} d^{2}(x, x_{\pm}),\\
\nonumber
\langle \nabla V(x), x-x_{\pm}\rangle &\ge \lambda^{2} d^{2}(x, x_{\pm}),\\
\langle H^{V}(x) v, v\rangle &\ge \lambda^{2} \vert v \vert^{2}.
\end{align}
Assuming the solution $q_{k}$ is constructed we introduce a time $\hat \xi_{k}$ defined as
\begin{equation}
\hat \xi_{k} = \min \left\{\xi\in (0, \xi_{k}): q_{k}(s)\in B_{R_{\lambda}}(x_{+})\,\, 
\forall\,\, s\in (\xi, \xi_{k}]\,\,\, and\,\,\,
q_{k}(s)\in B_{R_{\lambda}}(x_{-})\,\, 
\forall\,\, s\in [-\xi_{k}, -\xi)\right\}.
\end{equation}
We also represent 
$\Omega_{k} = \hat \Omega_{k}^{0} \cup \hat \Omega_{k}^{+} \cup \hat \Omega_{k}^{-}$ and 
$\Omega_{k+1} = \Omega_{k} \cup \Omega_{k+1}^{+} \cup \Omega_{k+1}^{-}$, where
\begin{align}
\nonumber
&\hat \Omega_{k}^{0} = [-\hat \xi_{k}, \hat \xi_{k}],\,\,\, 
\hat \Omega_{k}^{+} = [\hat \xi_{k}, \xi_{k}],\,\,\,
\hat \Omega_{k}^{-} = [- \xi_{k}, -\hat \xi_{k}],\\
&\Omega_{k+1}^{+} = [\xi_{k}, \xi_{k+1}],\,\,\,
\Omega_{k+1}^{-} = [- \xi_{k+1}, - \xi_{k}].
\end{align}

Consider the equations of motion (2.4) on the intervals $\hat\Omega_{k}^{\pm}$.  In the local chart $U^{\pm}$ one may rewrite them as
\begin{eqnarray}
%\nonumber
\frac{{\rm d}^{2}v^{\pm, j}}{{\rm d}\xi^{2}} + 
\Gamma^{j}_{m,l}(v^{\pm})\frac{{\rm d}v^{\pm, m}}{{\rm d}\xi}\frac{{\rm d}v^{\pm, l}}{{\rm d}\xi} + p(\xi)\frac{{\rm d}v^{\pm, j}}{{\rm d}\xi}+ 
\sigma(\xi) \left({\rm grad}V(v^{\pm})\right)^{j}=0,\,\,\,
j = 1,\ldots, n.
\end{eqnarray}
We consider only the case of the interval $\hat \Omega_{k}^{+}$ and skip the index $'+'$ for simplicity. The interval $\hat \Omega_{k}^{-}$ can be studied in a similar way. For any solution $v$ of the equations (5.9) define 
\begin{equation}
\rho(\xi) = \left(g_{m l}(v(\xi))v^{m}(\xi) v^{l}(\xi)\right)^{1/2}.
\end{equation}
Then, due to (5.8) and (5.6) one gets
\begin{align*}
%\nonumber
&\frac{{\rm d}}{{\rm d}\xi} \rho(\xi) = 
\rho^{-1}(\xi) g_{m l}(v(\xi)) \frac{{\rm d}v^{m}}{{\rm d}\xi} v^{ l}(\xi),
\\
%\nonumber
&\frac{1}{r(\xi)}\frac{{\rm d}}{{\rm d}\xi}r(\xi) \frac{{\rm d}}{{\rm d}\xi} \rho(\xi) = 
\rho^{-1}(\xi)\biggl(
g_{m l}(v(\xi)) \frac{{\rm d}v^{m}}{{\rm d}\xi}\frac{{\rm d}v^{ l}}{{\rm d}\xi} -
\rho^{-2}(\xi) \left( g_{m l}(v(\xi)) \frac{{\rm d}v^{m}}{{\rm d}\xi} v^{l}(\xi)\right)^{2} - 
\\
%\nonumber
&\hspace{5.6cm}
\sigma(\xi) g_{m l}(v(\xi)) \left({\rm grad}V(v(\xi))\right)^{m} v^{l}(\xi)\biggr)\ge \lambda^{2}\rho(\xi).
\end{align*}
If we take the solution $v_{k} = q_{k}\vert_{\hat\Omega_{k}}$ and substitute to (5.10) then $\rho$ will satisfy the following boundary conditions:
\begin{equation*}
\rho(\hat\xi_{k}) = R_{\lambda},\quad \rho(\xi_{k}) = 0.
\end{equation*}
Hence, by Collarary 1
\begin{equation*}
\bigl\vert \rho'(\xi_{k})\bigr\vert \le \frac{\lambda R_{\lambda}}{\sinh(\lambda(\xi_{k} - \hat\xi_{k}))}.
\end{equation*}
On the other hand, 
\begin{equation*}
\frac{v_{k}(\xi)}{\rho(\xi)}\to \frac{v'(\xi_{k})}{\rho'(\xi_{k})}\quad {\rm as}\,\, \xi\to \xi_{k}
\end{equation*}
implies
\begin{equation*}
\bigl\vert v'_{k}(\xi_{k})\bigr\vert = 
\left(g_{m l}(v_{k}(\xi_{k}))\frac{{\rm d}v^{m}_{k}}{{\rm d}\xi}(\xi_{k})
\frac{{\rm d}v^{ l}_{k}}{{\rm d}\xi}(\xi_{k})\right)^{1/2} = \bigl\vert \rho'(\xi_{k})\bigr\vert.
\end{equation*}
Thus, one finally has
\begin{equation}
\bigl\vert q'_{k}(\xi_{k})\bigr\vert \le \frac{\lambda R_{\lambda}}{\sinh(\lambda(\xi_{k} - \hat\xi_{k}))}.
\end{equation}
Note that the same estimate holds also for $\bigl\vert q'_{k}(-\xi_{k})\bigr\vert$.

This leads to the following
\begin{lemmas}
Let $q_{k}\in Q_{k}(\mathcal{M}, x_{-}, x_{+})$ be non-singular critical point of $I_{k}$ and $\varphi\in T_{q_{k}}\mathfrak{M}_{k+1}$ then 
\begin{equation}
\nonumber
\Bigl\vert I'_{k+1}[q_{k}](\varphi)\Bigr\vert \le 
b_{k} \Vert \varphi \Vert_{k+1},
\end{equation}
where 
\begin{equation}
b_{k} = \left(\bigl(r(\xi_{k})\bigr)^{1/2} + \bigl(r(-\xi_{k})\bigr)^{1/2}\right)
\frac{\lambda R_{\lambda}(\xi_{k+1}-\xi_{k})^{1/2}}{\sinh(\lambda(\xi_{k} - \hat\xi_{k}))}.
\end{equation}
\end{lemmas}
PROOF:- Note that $q_{k}\in C^{2}(\Omega_{k})$ and satisfies (2.4). Then using (3.2) and integrating by parts give
\begin{align}
\nonumber
I'_{k+1}[q_{k}](\varphi) = &\int\limits_{\Omega_{k}}\left(
\langle D_{\xi}q_{k}(\xi), D_{\xi}\varphi(\xi)\rangle - 
\sigma(\xi)\langle{\rm grad}V(q_{k}(\xi)), \varphi(\xi)\rangle\right)r(\xi) {\rm d}\xi +\\
\nonumber
&\int\limits_{\Omega_{k+1}\setminus \Omega_{k}}\left(
\langle D_{\xi}\chi(\xi), D_{\xi}\varphi(\xi)\rangle - 
\sigma(\xi)\langle{\rm grad}V(\chi(\xi)), \varphi(\xi)\rangle\right)r(\xi) {\rm d}\xi =\\
\nonumber
&-\int\limits_{\Omega_{k}}\langle D_{\xi}r(\xi)D_{\xi}q_{k}(\xi) + 
\sigma(\xi)r(\xi){\rm grad}V(q_{k}(\xi)), \varphi(\xi)\rangle {\rm d}\xi + 
r(\xi)\langle D_{\xi}q_{k}(\xi), \varphi(\xi)\rangle\biggl\vert_{-\xi_{k}}^{\xi_{k}} =\\ 
&r(\xi)\langle q'_{k}(\xi), \varphi(\xi)\rangle\biggl\vert_{-\xi_{k}}^{\xi_{k}}.
\end{align}
Applying the Schwartz inequality one gets
\begin{equation}
\vert \varphi(\xi_{k})\vert \le \int\limits_{\Omega_{k}^{+}}\vert \varphi'(\xi)\vert {\rm d}\xi \le \left(\int\limits_{\Omega_{k}^{+}}r^{-1}(\xi){\rm d}\xi\right)^{1/2}
\left(\int\limits_{\Omega_{k}^{+}}\vert \varphi'(\xi)\vert^{2}r(\xi){\rm d}\xi\right)^{1/2}\le
\left(\int\limits_{\Omega_{k}^{+}}r^{-1}(\xi){\rm d}\xi\right)^{1/2} \Vert \varphi\Vert_{k+1}.
\end{equation}
Hence, (5.13) together with (5.11), (5.12) yield
\begin{align*}
%\nonumber
\Bigl\vert I'_{k+1}[q_{k}](\varphi) \Bigr\vert \le &
r(\xi_{k})\bigl\vert q'_{k}(\xi_{k})\bigr\vert \cdot \bigl\vert \varphi(\xi_{k})\bigr\vert \le\\
&\left(r(\xi_{k})\left(\int\limits_{\Omega_{k}^{+}}r^{-1}(\xi){\rm d}\xi\right)^{1/2} + 
r(-\xi_{k})\left(\int\limits_{\Omega_{k}^{-}}r^{-1}(\xi){\rm d}\xi\right)^{1/2}\right)
\frac{\lambda R_{\lambda}}{\sinh(\lambda(\xi_{k} - \hat\xi_{k}))} \Vert \varphi \Vert_{k+1}.
\end{align*}
We finishes the proof by noting that $r$ is an increasing function. $\square$

To estimate the norm of the second derivative $I''_{k+1}[q_{k}]$ we study the variational equations along $q_{k}$ on the intervals $\hat\Omega_{k}^{+}(b)=[\hat\xi_{k}, b]$ and $\hat\Omega_{k}^{-}(b)=[-b, -\hat\xi_{k}]$ with $b\in [\xi_{k}, \xi_{k+1}]$. As in the previous case we consider only the interval $\hat\Omega_{k}^{+}(b)$ and skip the index $'+'$. 
Let $v = q_{k}\vert_{\hat\Omega_{k}^{+}(b)}$ be expressed in local coordinates of the chart $U^{+}$. Then the variational equations take the form
\begin{eqnarray}
\nonumber
\frac{{\rm d}^{2}w^{ j}}{{\rm d}\xi^{2}} + 
\Gamma^{j}_{m,l}(v)\frac{{\rm d}v^{m}}{{\rm d}\xi}\frac{{\rm d}w^{l}}{{\rm d}\xi} + p(\xi)\frac{{\rm d}w^{j}}{{\rm d}\xi} + 
R^{j}_{m,l,i}(v) \frac{{\rm d}v^{m}}{{\rm d}\xi}\frac{{\rm d}v^{l}}{{\rm d}\xi}w^{i} +
\sigma(\xi) \left(H^{V}(v)\right)^{j}_{m}w^{m}=0,\,\,\,
j = 1,\ldots, n
\end{eqnarray}
or equivalently
\begin{eqnarray}
%\nonumber
\frac{{\rm d}^{2}w^{ j}}{{\rm d}\xi^{2}} + 
\Gamma^{j}_{m,l}(v)\frac{{\rm d}v^{m}}{{\rm d}\xi}\frac{{\rm d}w^{l}}{{\rm d}\xi} + p(\xi)\frac{{\rm d}w^{j}}{{\rm d}\xi} - 
\mathcal{A}^{j}_{m}(v) w^{m} = 0,\,\,\,
j = 1,\ldots, n
\end{eqnarray}
where 
\begin{equation}
\mathcal{A}^{j}_{m}(v) = 
-R^{j}_{m,l,i}(v) \frac{{\rm d}v^{l}}{{\rm d}\xi}\frac{{\rm d}v^{i}}{{\rm d}\xi} -
\sigma(\xi) \left(H^{V}(v)\right)^{j}_{m}.
\end{equation}

We assume the operator $\mathcal{A}$ satisfies for some $\mu > 0$
\begin{equation}
\langle \mathcal{A} w, w\rangle \ge \mu^{2} \vert w\vert^{2}\quad 
\forall\,\, w\in T_{v(\xi)}\mathcal{M}, \,\,\,\forall\,\, \xi\in \hat\Omega_{k}^{\pm}.
\end{equation}
Although the second term in the r.h.s. of (5.16) is positive due to (5.6), one may note that the assumption (5.17) is rather restrictive. Since
\begin{equation*}
R^{j}_{m,l,i}(v) \frac{{\rm d}v^{m}}{{\rm d}\xi}\frac{{\rm d}v^{l}}{{\rm d}\xi}w^{m}w^{j} =
K(v',w)\left(\vert v'\vert^{2} \vert w\vert^{2} - \langle v', w\rangle^{2}\right), 
\end{equation*}
where $K$ denotes the sectional curvature, we conclude that (5.17) holds if for example the sectional curvature is non-positive in the ball $B_{R_{\lambda}}(x_{\pm})$. Another condition which yields (5.17) is the following. Consider 
\begin{equation*}
E_{k}(\xi) = \frac{1}{2}\vert q'_{k}\vert^{2} + \sigma(\xi)\left(V(q_{k}(x)) - V(\chi(\xi)) \right).
\end{equation*} 
Since $q_{k}$ solves (2.4)
\begin{equation*}
E'_{k}(\xi) = -p(\xi)\vert q'_{k}(\xi)\vert^{2}
\end{equation*}
and $E_{k}$ decreases on $\hat\Omega_{k}^{+}$. Besides, if $x\in B_{R_{\lambda}}(x_{+})$ $V(x_{+})-V(x) \le \frac{1}{2}\nu^{2} d^{2}(x, x_{+})$ for some $\nu > \max_{k}\{\Lambda_{k}^{+}\}$. Taking this into account one gets for any 
$\xi\in \hat\Omega_{k}^{\pm}$
\begin{equation*}
\vert q'_{k}(\xi)\vert^{2}\le 2 E_{k}(\hat \xi_{k}) + \nu^{2}R^{2}_{\lambda} \le
\vert q'_{k}(\hat\xi_{k})\vert^{2} + (\nu^{2} - \lambda^{2})R^{2}_{\lambda} \le
\left(\nu^{2}\coth^{2}\bigl(\nu(\xi_{k}-\hat\xi_{k})\bigr) +  
(\nu^{2} - \lambda^{2})\right)R^{2}_{\lambda}.
\end{equation*}
Introduce
\begin{equation*}
K_{\max} = \max\{K(w_{1}, w_{2}): w_{j}\in T_{x}\mathcal{M}, \,\,
x\in B_{R_{\lambda}}(x_{+}),\,\,  j=1,2\}.
\end{equation*}
Hence, the condition (5.17)  holds if
\begin{equation}
K_{\max} \left(\nu^{2}\coth^{2}\bigl(\nu(\xi_{k}-\hat\xi_{k})\bigr) +  
(\nu^{2} - \lambda^{2})\right)R^{2}_{\lambda} \le \lambda^{2} - \mu^{2}.
\end{equation}

We consider the variational equations (5.15) and supply them by boundary conditions
\begin{equation}
w(\hat\xi_{k}) = w_{0},\quad w(b) = 0,\quad \vert w_{0}\vert = 1.
\end{equation}
Denote by $\hat w = \hat w_{b}(\xi, b)$ the solution of the boundary value problem (5.15), (5.19). One may note that due to (5.17) $\hat w$ depends smoothly on the parameter $b$ and gives the minimum to the functional
\begin{equation*}
\hat I_{k, b}[w] = \int\limits_{\hat\xi_{k}}^{b}\biggl(\langle D_{\xi}w(\xi), D_{\xi}w(\xi)\rangle + 
\langle \mathcal{A}w(\xi), w(\xi)\rangle\biggr)r(\xi){\rm d}\xi
\end{equation*}
defined on the set of absolutely continuous functions $AC(\hat\Omega_{k}(b), T\mathcal{M})$ satisfying the boundary conditions (5.19). Introduce
\begin{equation}
\hat I_{k}(b) = \hat I_{k,b}[\hat w].
\end{equation}
Differentiating (5.20) with respect to $b$ and taking into account that $\hat w$ solves (5.15) one gets
\begin{eqnarray}
\nonumber
\frac{{\rm d} \hat I_{k}(b)}{{\rm d}b} = \vert \partial_{\xi} \hat w(b,b)\vert^{2} r(b) +
2\int\limits_{\hat \Omega_{k}(b)}
\left(\langle D_{\xi}\hat w(\xi, b),\partial_{b} D_{\xi} \hat w(\xi, b)\rangle +
\langle \mathcal{A} \hat w(\xi, b), \partial_{b} \hat w(\xi, b)\rangle\right)r(\xi){\rm d}\xi = \\
%\nonumber
\vert \partial_{\xi} \hat w(b,b)\vert^{2} r(b) + 
2\langle D_{\xi}\hat w(\xi, b), \partial_{b} \hat w(\xi, b)\rangle\biggl\vert_{\hat \xi_{k}}^{b} =
-\vert \partial_{\xi} \hat w(b,b)\vert^{2} r(b).
\end{eqnarray}
In the latter equality we used
\begin{equation*}
\partial_{b} \hat w(0,b) = 0,\quad \partial_{\xi} \hat w(b,b) + \partial_{b} \hat w(b,b) = 0
\end{equation*}
what follows from (5.19).
\begin{lemmas}
Let $q_{k}\in Q_{k}(\mathcal{M}, x_{-}, x_{+})$ be a non-singular critical point of $I_{k}$ such that 
$$
I''_{k}[q_{k}](\psi, \psi)\ge \gamma_{k}\Vert \psi\Vert_{k}^{2}\quad 
\forall\,\, \psi\in T_{q_{k}}\mathfrak{M}_{k}.
$$
Then for any $\varphi\in T_{q_{k}}\mathfrak{M}_{k+1}$  
\begin{equation}
I''_{k+1}[q_{k}](\varphi, \varphi)\ge a_{k}\Vert \varphi\Vert_{k+1}^{2},\,\,\,
a_{k} = \gamma_{k} - \Delta_{k}
\end{equation}
where
\begin{eqnarray}
\Delta_{k} = \frac{(1-\gamma_{k})\mu^{2}}{\sinh^{2}\bigl(\mu(\xi_{k}-\hat\xi_{k})\bigr)}
\left(
\int\limits_{\hat\Omega_{k}^{+}\cup \Omega_{k+1}^{+}}r^{-1}(\xi){\rm d}\xi \,
\int\limits_{\Omega_{k+1}^{+}}r(\xi){\rm d}\xi +
\int\limits_{\hat\Omega_{k}^{-}\cup \Omega_{k+1}^{-}}r^{-1}(\xi){\rm d}\xi \,
\int\limits_{\Omega_{k+1}^{-}}r(\xi){\rm d}\xi\right).
\end{eqnarray}
\end{lemmas}
PROOF:- For any $\varphi\in T_{q_{k}}\mathfrak{M}_{k+1}$ such that 
$\varphi(\pm \hat \xi_{k})\neq 0$ we consider a function $\hat\varphi\in T_{q_{k}}\mathfrak{M}_{k}$, which satisfies $\hat\varphi(\xi) = \varphi(\xi)$ for all
$\xi\in \hat\Omega_{k}^{0}$ (see (5.8) for the definition).
Then one gets
\begin{eqnarray}
\nonumber
I''_{k+1}[q_{k}](\varphi, \varphi) = I''_{k+1}[q_{k}](\hat\varphi, \hat\varphi) + 
I''_{k+1}[q_{k}](\varphi, \varphi) - I''_{k+1}[q_{k}](\hat\varphi, \hat\varphi) \ge\\
%\nonumber
\gamma_{k}\Vert \varphi\Vert_{k+1}^{2} + (1-\gamma_{k})\Bigl(
B_{k+1}[q_{k}](\varphi, \varphi) - B_{k+1}[q_{k}](\hat\varphi, \hat\varphi)\Bigr),
\end{eqnarray}
where the functional $B_{k+1}$ is defined by
\begin{eqnarray}
\nonumber
B_{k+1}[q](\varphi, \varphi) = 
\int\limits_{\Omega_{k+1}}\biggl(
\langle D_{\xi}\varphi(\xi), D_{\xi}\varphi(\xi)\rangle - (1-\gamma_{k})^{-1} \biggl[
\langle R(D_{\xi}q(\xi),\varphi(\xi))\varphi(\xi), D_{\xi}q(\xi)\rangle + \\
\nonumber
\langle\left(\sigma(\xi)\nabla V(q(\xi)) + \gamma_{k}I\right)\varphi(\xi), \varphi(\xi)\rangle \biggr]\biggr){\rm d}\xi = \\
\nonumber
\int\limits_{\Omega_{k+1}}\biggl(
\langle D_{\xi}\varphi(\xi), D_{\xi}\varphi(\xi)\rangle + 
\langle\mathcal{B}[q] \varphi(\xi), \varphi(\xi)\rangle \biggr)
{\rm d}\xi,
\end{eqnarray}
with
\begin{eqnarray}
\nonumber
\mathcal{B}[q] = (1-\gamma_{k})^{-1}\left(\mathcal{A}[q] - \gamma_{k}I\right)
\end{eqnarray}
Due to assumption (5.17) the following estimate holds
\begin{eqnarray}
\nonumber
\langle\mathcal{B}[q_{k}]\varphi(\xi), \varphi(\xi)\rangle\ge
\frac{\mu^{2}-\gamma_{k}}{1-\gamma_{k}}\vert \varphi(\xi)\vert^{2}\ge
\mu^{2} \vert\varphi(\xi)\vert^{2},\quad \forall\,\, \xi\in \hat\Omega_{k}^{\pm}\cup\Omega_{k+1}^{\pm}.
\end{eqnarray}
We consider the equations (5.15) with the operator $\mathcal{A}$ replaced by $\mathcal{B}$ and supply them by the boundary conditions (5.19) on the interval $\hat\Omega_{k}^{\pm}(\xi_{k})$ with parameter $w_{0}$ chosen as
\begin{equation*}
w_{0} = \varphi(\pm\hat \xi_{k})/\vert \varphi(\pm\hat \xi_{k})\vert.
\end{equation*}
Denote the corresponding solution by $\hat w^{\pm}$. It minimizes the functional $B_{k+1}$ restricted to the set of absolutely continuous functions $AC(\hat\Omega_{k}^{\pm}(\xi_{k}), T\mathcal{M})$ satisfying the boundary conditions (5.19). Using the solution $\hat w^{\pm}$, define 
$\hat\varphi\in T_{q_{k}}\mathfrak{M}_{k}$ as
\begin{equation*}
\hat \varphi(\xi)=
\begin{cases} 
\vert \varphi(-\hat \xi_{k})\vert \cdot\hat w^{-}(\xi),\,\, \xi\in \hat\Omega_{k}^{-},\\
\varphi(\xi),\,\, \xi\in \hat \Omega_{k}^{0},\\
\vert \varphi(\hat \xi_{k})\vert \cdot \hat w^{+}(\xi),\,\, \xi\in \hat\Omega_{k}^{+}.
\end{cases}
\end{equation*}
Note here that in the case $\varphi(\pm\hat \xi_{k}) = 0$ we continue $\hat\varphi$ by $0$ on the interval $\hat\Omega_{k}^{\pm}$. 
Then taking into account (5.21) and Collorary 1 to Lemma 6 one obtains
\begin{align}
\nonumber
\left[ B_{k+1}[q_{k}](\varphi, \varphi) - B_{k+1}[q_{k}](\hat\varphi, \hat\varphi)\right]_{-}\le& \\
\frac{\mu^{2}R_{\mu}^{2}}{\sinh^{2}\bigl(\mu(\xi_{k}-\hat\xi_{k})\bigr)}
&\left(\vert \varphi(\hat \xi_{k})\vert^{2}\int\limits_{\Omega_{k+1}^{+}}r(\xi){\rm d}\xi +
\vert \varphi(-\hat \xi_{k})\vert^{2}\int\limits_{\Omega_{k+1}^{+}}r(\xi){\rm d}\xi\right),
\end{align}
where $[x ]_{-}$ stands for the negative part of a real number $x$.

Applying the same arguments as in (5.14) yields
\begin{equation}
\vert \varphi(\pm \hat\xi_{k})\vert^{2} \le \int\limits_{\hat\Omega_{k}^{\pm}\cup \Omega_{k+1}^{\pm}}r^{-1}(\xi){\rm d}\xi \cdot \Vert\varphi \Vert_{k+1}^{2}.
\end{equation}
Finally substituting (5.25), (5.26) into (5.24)  finishes the proof. $\square$

The next lemma provides information on the time $\hat \xi_{k}$. Let $\hat{\mathfrak{M}}_{k}$ be the set of absolutely continuous functions 
$\hat{\mathfrak{M}}_{k} = \{q\in AC([-1/2, 1/2], \mathcal{M}): q(\pm 1/2)\in S_{R_{\lambda}}(x_{\pm})\}$ where $S_{R}(x)$ denotes the sphere of radius $R$ centered at $x\in \mathcal{M}$. We define 
\begin{equation*}
L_{k}^{2} = \min\limits_{q\in \hat{\mathfrak{M}}_{k}}\int\limits_{-1/2}^{1/2}\vert \dot q(s)\vert^{2}{\rm d}s.
\end{equation*}
We also introduce $\hat I_{k} = I_{k}[q_{k}]$ and 
\begin{equation}
\eta_{k} = \int\limits_{\Omega_{k}}r^{-1}(\xi){\rm d}\xi,\quad
h_{k} = \frac{2 R_{\lambda}^{2}\eta_{k}}{2\hat I_{k}\eta_{k} - L_{k}^{2}  + 
\left(\left(2\hat I_{k}\eta_{k} - L_{k}^{2}\right)^{2} - 
8 R_{\lambda}^{2}\hat I_{k}\eta_{k}\right)^{1/2}}.
\end{equation}
\begin{lemmas}
Let $q_{k}\in Q_{k}(\mathcal{M}, x_{-}, x_{+})$ be a critical point of $I_{k}$. Suppose there exists a positive solution $\zeta_{k}^{\pm}$ of the equation
$$
\pm\int\limits_{\pm\zeta_{k}^{\pm}}^{\pm\xi_{k}}r^{-1}(\xi){\rm d}\xi = h_{k}.
$$
Then $\hat\xi_{k} \le \zeta_{k}$, where $\zeta_{k} = \max\{\zeta_{k}^{\pm}\}$.
\end{lemmas}
PROOF:- First we observe that
\begin{eqnarray}
\nonumber
\int\limits_{\xi}^{\xi_{k}}\vert q'_{k}(s)\vert^{2}r(s){\rm d}s = 
\int\limits_{\Omega_{k}}\vert q'_{k}(s)\vert^{2}r(s){\rm d}s - 
\int\limits_{-\xi_{k}}^{\xi}\vert q'_{k}(s)\vert^{2}r(s){\rm d}s \le\\
\nonumber
2 \hat I_{k} - \min\limits_{q\in \hat{\mathfrak{M}}_{k}^{+}(\xi)}
\int\limits_{-\xi_{k}}^{\xi}\vert q'(s)\vert^{2}r(s){\rm d}s\le
2 \hat I_{k} - \left(\int\limits_{-\xi_{k}}^{\xi}r^{-1}(s){\rm d}s\right)^{-1}
\Bigl(L_{k}^{2} + R_{\lambda}^{2}\Bigr),
\end{eqnarray}
where
\begin{equation*}
\hat{\mathfrak{M}}_{k}^{+}(\xi) = \{q\in AC([-\xi_{k}, \xi], \mathcal{M}): q(-\xi_{k})=x_{-},\,\, q(\xi)\in S_{R_{\lambda}}(x_{+})\}.
\end{equation*}

Hence,
\begin{equation}
d^{2}(q_{k}(\xi), x_{+})\le \int\limits_{\xi}^{\xi_{k}}r^{-1}(s){\rm d}s
\left(2\hat I_{k} - \left(\int\limits_{-\xi_{k}}^{\xi}r^{-1}(s){\rm d}s\right)^{-1}
\Bigl(L_{k}^{2} + R_{\lambda}^{2}\Bigr)\right).
\end{equation}
In a similar way one may show that
\begin{equation}
d^{2}(q_{k}(\xi), x_{-})\le \int\limits_{-\xi_{k}}^{\xi}r^{-1}(s){\rm d}s
\left(2\hat I_{k} - \left(\int\limits_{\xi}^{\xi_{k}}r^{-1}(s){\rm d}s\right)^{-1}
\Bigl(L_{k}^{2} + R_{\lambda}^{2}\Bigr)\right).
\end{equation}

For $\zeta\ge 0$ we consider the following equation
\begin{equation*}
\pm\int\limits_{\pm\zeta}^{\pm\xi_{k}}r^{-1}(s){\rm d}s
\left(2\hat I_{k} - \left(\pm\int\limits_{\mp\xi_{k}}^{\pm\zeta}r^{-1}(s){\rm d}s\right)^{-1}
\Bigl(L_{k}^{2} + R_{\lambda}^{2}\Bigr)\right) = 
R_{\lambda}^{2}. 
\end{equation*}
Introducing  
\begin{equation*}
h = \pm\int\limits_{\pm\zeta}^{\pm\xi_{k}}r^{-1}(s){\rm d}s
\end{equation*}
one may rewrite this equation as a quadratic one with respect to $h$:
\begin{equation} 
2\hat I_{k} h^{2} - \Bigl(2\hat I_{k}\eta_{k} - L_{k}^{2}\Bigr)h + R_{\lambda}^{2}\eta_{k} = 0.
\end{equation}
Note that $h_{k}$ is the smallest solution of (5.30). Taking into account (5.28), (5.29) and definitions of $\zeta_{k}^{\pm}$, we conclude that $\zeta_{k}$ provides the upper bound for the time $\hat \xi_{k}$. $\square$

One may remark here that if there does not exist a positive solution $\zeta_{k}^{+}$ (resp. $\zeta_{k}^{-}$) then the inequality (5.28) (resp. (5.29)) holds for any positive (resp. negative) $\xi$ such that $\vert\xi\vert \le \xi_{k}$. In this case one may set $\zeta_{k}^{+}=0$  (resp. $\zeta_{k}^{-}=0$) to preserve the statment of the Lemma 9 without changes.

Finally, one may apply (5.21) and (5.11) to show that
\begin{equation}
\hat I_{k} < \hat I_{k-1} - \Delta\hat I_{k-1},\quad
\Delta\hat I_{k-1} = 
\frac{R_{\lambda}^{2} \nu^{2}}{2\sinh^{2}\Bigl(\nu(\xi_{k-1}-\hat\xi_{k-1})\Bigr)}
\int\limits_{\xi_{k-1}}^{\xi_{k}}r(\xi){\rm d}\xi,
\end{equation}
where $\nu > \max\{\Lambda_{j}^{\pm}, j=1,\ldots,n\}$ such that 
$V(x_{\pm}) - V(x)\le \frac{1}{2}\nu^{2}d^{2}(x, x_{\pm})$ for any 
$x\in B_{R_{\lambda}}(x_{\pm})$.
If one substitutes (5.31) into the definition (5.27), the value of $h_{k}$ decreases, while $\zeta_{k}$ increases. Hence, the following corollary takes place
\begin{corollaries}
Let $q_{k}\in Q_{k}(\mathcal{M}, x_{-}, x_{+})$ be a critical point of $I_{k}$. Define $\hat\zeta_{k}^{\pm}$ to be a positive solution of the equation
\begin{equation}
\pm\int\limits_{\pm\zeta_{k}^{\pm}}^{\pm\xi_{k}}r^{-1}(\xi){\rm d}\xi =
\frac{2 R_{\lambda}^{2}\eta_{k}}{2(\hat I_{k-1} - \Delta\hat I_{k-1})\eta_{k} - L_{k}^{2}  + 
\left(\left(2(\hat I_{k-1} - \Delta\hat I_{k-1})\eta_{k} - L_{k}^{2}\right)^{2} - 
8 R_{\lambda}^{2}\hat I_{k}\eta_{k}\right)^{1/2}}.
\end{equation}
if it exists or to be zero otherwise.
Then $\hat\xi_{k} \le \hat\zeta_{k}$, where $\hat\zeta_{k} = \max\{\hat\zeta_{k}^{\pm}\}$.
\end{corollaries}

To apply the Newton-Kantorovich theorem on  the $k$-th step one needs to check that the second deriavative of the functional $I_{k}$ is locally Lipschitz. First we note that there exist positive constants $C_{g}$, $C_{K}$ and $C_{V}$ such that for any $x\in \mathcal{M}$ and $v, w\in T_{x}\mathcal{M}$ 
\begin{eqnarray}
%\nonumber
\sum\limits_{i=1}^{n}\left\vert\Gamma^{i}_{jk}(x)v^{j}v^{k}\right\vert \le 
C_{g}\Vert v\Vert^{2},\quad
\Bigl\vert\langle R(v,w)w, v\rangle\Bigr\vert\le C_{K}\Vert v\Vert\cdot \Vert w\Vert,\quad
\left\vert\langle H^{V}(x)v, v\rangle\right\vert\le C_{V}\Vert v\Vert^{2}.
\end{eqnarray}

Assume $q\in \mathfrak{M}_{k}$ and consider an open ball $B(q, R)$ of radius $R$ centered at $q$. Then we arrive at the following proposition
\begin{propositions}
The covariant derivative $\mathcal{D} \hat X_{I_{k}}$ is Lipschitz in $B(q, R)$ with constant 
$C_{L}^{I_{k}} = C_{L}^{I_{k}}(q)$
\begin{eqnarray}
C_{L}^{I_{k}} = 2\left(1 + \left(C_{g}^{2}+C_{K}\right)C_{r}^{2}\left(J_{k}[q] + 
\frac{1}{2} R^{2}\right) + \frac{1}{4}C_{V}C_{r}^{2}\right),
\end{eqnarray}
where
$$
J_{k}[q] = \frac{1}{2}\int\limits_{\Omega_{k}}\vert q'(\xi)\vert^{2}r(\xi){\rm d}\xi.
$$
\end{propositions}
PROOF:-
For any $q_{*}\in B(q, R)$ we take a geodesic $\alpha: [0, 1]\to \mathfrak{M}_{k}$ such that $\alpha(0) = q$ and $\alpha(1) = q_{*}$. Besides we consider an arbitrary $\alpha$-parallel vector field $X_{\varphi}\in \mathfrak{X}(\mathfrak{M}_{k})$ along the curve $\alpha$, i.e. $X_{\varphi}(\alpha(s)) = \hat P_{\alpha, 0, s} X_{\varphi}(\alpha(0))$. Then the Leibniz's formula reads (see e.g. \cite{FerSva}):
\begin{equation*}
\hat P_{\alpha, s, 0}(X_{\varphi}(\alpha(s)) = X_{\varphi}(\alpha(s)) + 
\int\limits_{0}^{s}\hat P_{\alpha, p, 0}\left( \mathcal{D}X_{\varphi}(\alpha(p))\alpha'(p)\right){\rm d}p.
\end{equation*}
Since $X_{\varphi}$ is $\alpha$-parallel and
$\mathcal{D}\hat X_{I_{k}}(\alpha(s))X_{\varphi}(\alpha(s)) = 
(\nabla_{X_{\varphi}}\hat X_{I_{k}})(\alpha(s))$ one gets
\begin{equation*}
\hat P_{\alpha, 1, 0} \mathcal{D} \hat X_{I_{k}}(\alpha(1)) \hat P_{\alpha, 0, 1} X_{\varphi}(\alpha(0)) - 
\mathcal{D} \hat X_{I_{k}}(\alpha(0)) X_{\varphi}(\alpha(0)) =
\int\limits_{0}^{1}\hat P_{\alpha, s, 0}\left(\nabla_{\alpha'(s)} \left(\nabla_{X_{\varphi}}\hat X_{I_{k}}(\alpha(s))\right)\right)
{\rm d}s.
\end{equation*}
Let $\psi\in T_{q}\mathfrak{M}_{k}$ and $X_{\psi}(\alpha(s)) = \hat P_{\alpha, 0, s}\psi$. Then
\begin{eqnarray}
\nonumber
\Bigl\langle\Bigl\langle 
\int\limits_{0}^{1}\hat P_{\alpha, s, 0}\left(\nabla_{\alpha'(s)} \left(\nabla_{X_{\varphi}}\hat X_{i_{k}}(\alpha(s))\right)\right)
{\rm d}s, \psi\Bigr\rangle\Bigr\rangle =
\int\limits_{0}^{1} \Bigl\langle\Bigl\langle 
\nabla_{\alpha'(s)} \left(\nabla_{X_{\varphi}}\hat X_{I_{k}}(\alpha(s))\right), X_{\psi}(\alpha(s))
\Bigr\rangle\Bigr\rangle {\rm d}s = \\
\nonumber
\int\limits_{0}^{1} X_{\alpha'(s)} \Bigl\langle\Bigl\langle 
\nabla_{X_{\varphi}}\hat X_{I_{k}}(\alpha(s)), X_{\psi}(\alpha(s)) \Bigr\rangle\Bigr\rangle 
{\rm d}s - 
\int\limits_{0}^{1} \Bigl\langle\Bigl\langle 
\nabla_{X_{\varphi}}\hat X_{I_{k}}(\alpha(s)), \nabla_{\alpha'(s)} X_{\psi}(\alpha(s)) \Bigr\rangle\Bigr\rangle {\rm d}s =\\
\nonumber
\int\limits_{0}^{1} X_{\alpha'(s)} \Bigl\langle\Bigl\langle 
\mathcal{D} \hat X_{I_{k}}(\alpha(s)) X_{\varphi}(\alpha(s)), X_{\psi}(\alpha(s)) \Bigr\rangle\Bigr\rangle {\rm d}s.
\end{eqnarray}
Hence
\begin{align}
\nonumber
\left\vert
\Bigl\langle\Bigl\langle 
\int\limits_{0}^{1}\hat P_{\alpha, s, 0}\left(\nabla_{\alpha'(s)} \left(\nabla_{X_{\varphi}}\hat X_{I_{k}}(\alpha(s))\right)\right)
{\rm d}s, \psi\Bigr\rangle\Bigr\rangle
\right\vert \le & \\
\int\limits_{0}^{1} \Vert \alpha'(s)\Vert \cdot &
\left\vert \Bigl\langle\Bigl\langle 
\mathcal{D}\hat X_{I_{k}}(\alpha(s)) X_{\varphi}(\alpha(s)), X_{\psi}(\alpha(s))\Bigr\rangle\Bigr\rangle \right\vert {\rm d}s.
\end{align}
Note that for any $\varphi, \psi \in T_{q}\mathfrak{M}_{k}$
\begin{equation*}
\langle\langle 
\mathcal{D}\hat X_{I_{k}}(q) X_{\varphi}, X_{\psi}
\rangle\rangle =
X_{\varphi}  \langle\langle \hat X_{I_{k}}(q) , X_{\psi} \rangle\rangle - 
\langle\langle \hat X_{I_{k}}(q), \hat\nabla_{\varphi} X_{\psi}
\rangle\rangle = 
X_{\varphi} \left(I'_{k}[q](\psi) \right) - I'_{k}[q]\left( \hat\nabla_{\varphi} X_{\psi} \right).
\end{equation*}

Taking this into account (5.33) together with (2.1), (3.3)  and applying the Schwartz inequality one obtains
\begin{equation}
\left\vert\langle\langle 
\mathcal{D}\hat X_{I_{k}}(q) X_{\varphi}, X_{\psi}
\rangle\rangle\right\vert\le
2\left(1 + \left(C_{g}^{2} + C_{K}\right)C_{r}^{2}J_{k}[q] + 
\frac{1}{4}C_{V}C_{r}^{2}\right)\Vert X_{\varphi}\Vert\cdot \Vert X_{\psi}\Vert.
\end{equation}

Substitute (5.36) into (5.35) to obtain
\begin{align*}
\nonumber
\left\vert
\Bigl\langle\Bigl\langle 
\int\limits_{0}^{1}\hat P_{\alpha, s, 0}\left(\nabla_{\alpha'(s)} \left(\nabla_{X_{\varphi}}\hat X_{I_{k}}(\alpha(s))\right)\right)
{\rm d}s, \psi\Bigr\rangle\Bigr\rangle
\right\vert \le &\\
2\int\limits_{0}^{1} \Bigl(1 + \left(C_{g}^{2}+C_{K}\right)C_{r}^{2}J_{k}[\alpha(s)] + &
\frac{1}{4}C_{V}C_{r}^{2}\Bigr) 
\Vert \alpha'(s)\Vert \cdot
\Vert X_{\varphi}(\alpha(s))\Vert \cdot \Vert X_{\psi}(\alpha(s))\Vert {\rm d}s.
\end{align*}
Taking into account that $\hat P_{\alpha, s, p}$ is an isometry and $\alpha(s)\in B(q, R)$ yields
\begin{equation*}
\Bigl\Vert \hat P_{\alpha, 1, 0} \mathcal{D} X_{I_{k}}(p) \hat P_{\alpha, 0, 1} - 
\mathcal{D} X_{I_{k}}(q)\Bigr\Vert_{op}\le
2\left(1 + \left(C_{g}^{2}+C_{K}\right)C_{r}^{2}\left(J_{k}[q] + 
\frac{1}{2} R^{2}\right) + 
\frac{1}{4}C_{V}C_{r}^{2}\right)
\int\limits_{0}^{1} \Vert \alpha'(s)\Vert {\rm d}s,
\end{equation*}
what finishes the proof. $\square$

Applying the Newton-Kantorovich theorem we arrive at the following lemma
\begin{lemmas}
Let $q_{k}\in Q_{k}(\mathcal{M}, x_{-}, x_{+})$ be a non-degenerate critical point of the functional $I_{k}$, which satisfies the conditions of the Lemmae 7-9. If $\xi_{k+1} > \xi_{k}$ is chosen such that
\begin{equation}
p_{k} = a_{k}^{2} b_{k} C_{L} < \frac{1}{2},
\end{equation}
where $a_{k}, b_{k}$ are defined by (5.12), (5.22),
then there exists a non-degenerate critical point $q_{k+1}$ of the functional $I_{k+1}$ such that 
\begin{equation*}
I''_{k+1}[q_{k+1}](\varphi, \varphi) \ge \gamma_{k+1}\Vert \varphi \Vert_{k+1}^{2},\quad
\forall\,\, \varphi\in T_{q_{k+1}}\mathfrak{M}_{k+1},
\end{equation*}
with
\begin{equation}
\gamma_{k+1} = \sqrt{1-2 p_{k}} a_{k}^{-1}.
\end{equation}
\end{lemmas}
Applying recurrently (5.38) and using (5.22), (5.23) one gets
\begin{equation}
\gamma_{k+1} = \sqrt{1-2 p_{k}}\bigl(\gamma_{k}-\Delta_{k}\bigr) =
\prod\limits_{j=0}^{k}\sqrt{1-2 p_{j}}\gamma_{0} - 
\sum\limits_{j=0}^{k}\Delta_{j}\prod\limits_{i=j}^{k}\sqrt{1-2 p_{i}}.
\end{equation}
If we set
\begin{equation*}
A_{0} = 1,\quad 
A_{k} = \prod\limits_{j=0}^{k-1} \frac{1}{\sqrt{1-2 p_{j}}}
\end{equation*}
one may rewrite (5.39) as
\begin{equation}
\gamma_{k+1} = 
A_{k+1}^{-1}\left(\gamma_{0} - \sum\limits_{j=0}^{k}A_{j}\Delta_{j}\right).
\end{equation}
We represent the Kantorovich's condition (5.37) as
\begin{align}
\nonumber
b_{k} = \frac{p_{k}}{2C_{L}A_{k}^{2}}\left(\gamma_{0} - \sum\limits_{j=0}^{k}A_{j}\Delta_{j}\right)^{2}=&\\
\frac{p_{k}}{2C_{L}A_{k}^{2}} & \left[\left(\gamma_{0} - \sum\limits_{j=0}^{k-1}A_{j}\Delta_{j}\right)^{2} -
2\left(\gamma_{0} - \sum\limits_{j=0}^{k-1}A_{j}\Delta_{j}\right)A_{k}\Delta_{k} + 
A_{k}^{2}\Delta_{k}^{2}\right] 
\end{align}
and note that only last two terms in the r.h.s. of (5.67) depend on $\xi_{k+1}$. 

Thus, taking into account definitions (5.12), (5.22), (5.23) one may consider the condition (5.41) as a definition of $\xi_{k+1}$ while the condition (5.32) as a definition of $\hat \xi_{k+1}$. 

To construct the solution $q_{0}$ we take a sufficiently small $\xi_{0}>0$ and a geodesic $\varGamma$ connecting $x_{-}$ and $x_{+}$ on the interval $\Omega_{0}$. We assume the geodesic $\varGamma$ to be non-degenerate and to satisfy
\begin{equation}
J''_{0}[\varGamma](\varphi, \varphi) \ge C_{\varGamma} 
\int\limits_{\Omega_{0}}\vert \varphi'(\xi)\vert^{2}r(\xi){\rm d}\xi \quad
\forall\,\, \varphi\in T_{\varGamma}\mathfrak{M}_{0}.
\end{equation}
The Schwartz inequality yields
\begin{equation*}
\int\limits_{\Omega_{0}}\vert \varphi(\xi)\vert^{2}r(\xi){\rm d}\xi \le
C_{\xi_{0}} \int\limits_{\Omega_{0}}\vert \varphi'(\xi)\vert^{2}r(\xi){\rm d}\xi,\quad
C_{\xi_{0}} = \frac{1}{2}\int\limits_{\Omega_{0}}r(\xi){\rm d}\xi
\int\limits_{\Omega_{0}}r^{-1}(\xi){\rm d}\xi.
\end{equation*}
Hence, we can take $\xi_{0}$ such that for any $\varphi\in T_{\varGamma}\mathfrak{M}_{0}$
\begin{align}
\nonumber
&\biggl\vert I'_{0}[\varGamma](\varphi)\biggr\vert \le b_{0}\Vert \varphi\Vert_{0},\quad
b_{0} = C_{D,1}^{V} C_{\xi_{0}},\\
\nonumber
&I''_{0}[\varGamma](\varphi, \varphi)\ge a_{0}\Vert \varphi\Vert_{0}^{2},\quad
a_{0} = \frac{C_{\varGamma} - C_{D,2}^{V}C_{\xi_{0}}}{1+C_{\xi_{0}}},\\
&p_{0} = a_{0}^{2}b_{0}C_{L} < \frac{1}{2},
\end{align}
where 
$$
C_{D,1}^{V} = \max\limits_{\xi\in \Omega_{0}}\Bigl\vert \nabla V(\varGamma(\xi))\Bigr\vert,\quad
C_{D,2}^{V} = \max\limits_{\xi\in \Omega_{0}}\Vert H^{V}(\varGamma(\xi))\Vert.
$$
Then the Newton-Kantorovich theorem guarantees the existence of the solution $q_{0}\in \mathfrak{M}_{0}$ connecting $x_{-}$ and $x_{+}$ on the interval $\Omega_{0}$. Moreover, one has
\begin{equation*}
I''_{0}[q_{0}](\varphi, \varphi)\ge \gamma_{0}\Vert \varphi\Vert_{0}^{2},\quad \forall\,\,
\varphi\in T_{q_{0}}\mathfrak{M}_{0},\quad \hat d(\varGamma, q_{0}) < a_{0}b_{0}.
\end{equation*}
We also note that
\begin{equation*}
d(q_{0}(\xi), x_{+}) \le d(\varGamma(\xi), x_{+}) + d(q_{0}(\xi), \varGamma(\xi)) \le
d(\varGamma(\xi), x_{+}) + \int\limits_{\xi}^{\xi_{0}}r^{-1}(s){\rm d}s \cdot 
\hat d(\varGamma, q_{0}).
\end{equation*}
For sufficiently small $\sigma>0$ define
\begin{equation*}
\xi_{\varGamma} = \min\{\xi\in \Omega_{0}: d(\varGamma(s), x_{+})\le (1-\sigma)R_{\lambda},\quad
\forall\,\, s\in [\xi, \xi_{0}]\}
\end{equation*}
and  consider an inequality
\begin{equation*}
(1-\sigma)R_{\lambda} + a_{0}b_{0}\int\limits_{\xi}^{\xi_{0}}r^{-1}(s){\rm d}s \le R_{\lambda}.
\end{equation*}
We take $\xi_{0}$ to be sufficiently small such that (5.47) holds for all $\xi\in [\xi_{\varGamma}, \xi_{0}]$ and, hence,
\begin{equation}
\hat \xi_{0} \le \xi_{\varGamma}.
\end{equation}
Summarizing all the results we get
\begin{theorems}
Let $\xi_{0}$ be a positive real number and $\varGamma$ be a non-degenerate geodesic connecting $x_{-}$ and $x_{+}$ on the interval $\Omega_{0}$. It is assumed that $\varGamma$ satisfies the estimate (5.44) and $\xi_{0}$ is small enough the estimates (5.43) and (5.47) to be valid. Let $\{p_{k}\}_{k=1}^{\infty}$ be a sequence of positive real numbers such that 
$$
p_{k}<1/2, \quad \exists\,\, 
\lim\limits_{k\to \infty} \prod\limits_{j=0}^{\infty}\left(\sqrt{1-2p_{j}}\right)^{-1} = 
A_{\infty} < \infty.
$$ 
If sequences $\{\xi_{k}\}_{k=1}^{\infty}$, $\{\hat\xi_{k}\}_{k=1}^{\infty}$ generated by the equations (5.41) and (5.32), respectively, satisfy
\begin{equation*}
\lim\limits_{k\to \infty} \xi_{k} = \infty, \quad \exists\,\,
\lim\limits_{k\to \infty}\sum_{j=0}^{k}A_{j}\Delta_{j} < \gamma_{0},
\end{equation*}
then there exists a sequence of solutions $q_{k}$ connecting $x_{-}$ and $x_{+}$ on the interval $\Omega_{k}$ and this sequence converges to a transversal doubly asymptotic solution $q_{\infty}$.
\end{theorems}

\section{Case $f(t) = t^{m}$}
\setcounter{equation}{0}

In this section we consider a special case when the factor $f$ is of the form $f(t) = t^{m}, m\in \mathbb{N}$. As it was mentioned in the introduction, factors of such kind appear in the study of Lagrangian systems with turning points and provide the main example of this paper. 

In terms of the variable $\xi$ the factor $f$ corresponds to
\begin{equation}
r(\xi) = \left(\frac{m+2}{2}\xi\right)^{\frac{m}{m+2}},\quad
\sigma(\xi) = \left({\rm sign}\, \xi\right)^{m}, \quad
p(\xi) = \frac{m}{m+2} \xi^{-1}.
\end{equation}

For this case we give more explicit conditions which guarantee the existence of transversal doubly asymptotic trajectories. We begin with simplification of the equation (5.32) which defines the parameter $\hat\xi_{k}$. For positive $h$ and $\xi$ consider an equation
\begin{equation*}
\int\limits_{\zeta}^{\xi}r^{-1}(s){\rm d}s = h.
\end{equation*} 
If $\int_{0}^{\xi}r^{-1}(s){\rm d}s \ge h$ there exists a unique positive solution $\zeta(\xi; h)$.
Introduce
\begin{equation*}
g(\xi) = \int\limits_{0}^{\xi}r^{-1}(s){\rm d}s.
\end{equation*}
The conditions $A_{1} - A_{3}$ imply 
\begin{equation}
g(0)= 0,\quad
g'(\xi) = r^{-1}(\xi) > 0,\quad
g''(\xi) = -\frac{r'(\xi)}{r^{2}(\xi)} < 0,\quad
g'''(\xi) = 2\frac{(r'(\xi))^{2}}{r^{3}(\xi)} - \frac{r''(\xi)}{r^{2}(\xi)} > 0\quad
\forall\,\, \xi>0.
\end{equation}
Then the solution $\zeta$ can be expressed as 
\begin{equation*}
\zeta(\xi; h) = g^{-1}(g(\xi)-h).
\end{equation*}
Note that due to (6.2) $\Upsilon(\xi) = \xi - \zeta(\xi; h) \to +\infty$ as $\xi\to +\infty$.
Taking into account the particular form of the factor $r(\xi)$ one gets
\begin{equation*}
g(\xi) = \alpha \xi^{\frac{2}{m+2}},\quad
\Upsilon(\xi) = 
\xi\left(1 - \left(1-\alpha^{-1}h \xi^{-\frac{2}{m+2}}\right)^{\frac{m+2}{2}}\right),\quad
\alpha = \left(\frac{m+2}{2}\right)^{\frac{2}{m+2}}.
\end{equation*}
Since $1 - (1 - s)^{\beta} \ge s$ provided $\beta > 1$ and $s\in [0,1]$ we obtain
\begin{equation}
\Upsilon(\xi) \ge \alpha^{-1} h \xi^{\frac{m}{m+2}},\quad
\forall\,\, \xi\ge \left(\alpha^{-1} h\right)^{\frac{m+2}{2}}.
\end{equation}
Hence, the following lemma holds
\begin{lemmas}
Let $q_{k}\in Q_{k}(\mathcal{M}, x_{-}, x_{+})$  be a critical point of the functional $I_{k}$ and  $h_{k}$ defined by (5.27) satisfies $h_{k} \ge h_{*} >0$. Then 
\begin{equation}
\xi_{k} - \hat\xi_{k} \ge 
\begin{cases}
\alpha^{-1} h_{*} \xi_{k}^{\frac{m}{m+2}}, \quad 
\xi_{k}\ge \left(\alpha^{-1} h_{*}\right)^{\frac{m+2}{2}},\\
\xi_{k}, \quad \xi_{k}\le \left(\alpha^{-1} h_{*}\right)^{\frac{m+2}{2}}.
\end{cases}
\end{equation}
\end{lemmas}
PROOF:- The statement of the lemma follows from Lemma 9 and (6.3).$\square$
\newline
We remark that all $h_{k}$ are uniformly bounded by
\begin{equation}
h_{k} > \frac{R_{\lambda}^{2}}{2 I_{0}[\varGamma]},\quad \forall\,\, k\ge 0.
\end{equation}

Substitute (6.4) into formulae (5.12), (5.22), (5.23) and take into account (6.1) to obtain
\begin{align}
\nonumber
& b_{k}\le \frac{2 \lambda R_{\lambda} \alpha^{1/2}\xi_{k}^{\frac{m}{2(m+2)}}}
{\sinh\Bigl(\lambda\alpha^{-1} h_{*}\xi_{k}^{\frac{m}{m+2}}\Bigr)}
\left(\xi_{k+1} - \xi_{k}\right)^{1/2},
\\
& \Delta_{k} \le \frac{(1-\gamma_{k})\mu^{2}}
{\sinh^{2}\Bigl(\mu \alpha^{-1} h_{*}\xi_{k}^{\frac{m}{m+2}}\Bigr)}
\left(\frac{m+2}{2}\right)^{\frac{2(m+1)}{m+2}}\frac{1}{m+1}
\left( h_{*} + \alpha \left(\xi_{k+1}^{\frac{2}{m+2}} - \xi_{k}^{\frac{2}{m+2}}\right)\right)
\left(\xi_{k+1}^{\frac{2(m+1)}{m+2}} - \xi_{k}^{\frac{2(m+1)}{m+2}}\right).
\end{align}

One may note that the Newton-Kantorovich condition holds for any 
$\xi_{k+1}\in [\xi_{k}, \xi_{k+1}^{*}]$, where $\xi_{k+1}^{*}$ is defined be (5.37). Hence, without loss of generality we may additionally assume that 
\begin{equation*}
%\xi_{k+1} \le 2\xi_{k},\quad \forall\,\, k\ge 0.
\xi_{k+1} \le \xi_{k} + \xi_{k}^{\frac{m}{m+2}},\quad \forall\,\, k\ge 0.
\end{equation*}
This assumption together with inequalities:
\begin{equation*}
%\nonumber
(1+s)^{1+\beta} - 1 \le (2^{1+\beta}-1)s,\quad 
(1+s)^{\beta} - 1 \le s,\quad
0\le \beta\le 1, s\in [0, 1]
\end{equation*}
leads to
\begin{align}
\nonumber
& \Delta_{k}\le 
\frac{C_{\Delta}(1-\gamma_{k})\mu\alpha^{-1}h_{*}\xi_{k}^{\frac{m}{m+2}}}
{\sinh^{2}\Bigl(\mu \alpha^{-1} h_{*}\xi_{k}^{\frac{m}{m+2}}\Bigr)}
\left(\xi_{k+1} - \xi_{k}\right),\\
& C_{\Delta} = \left(1 + h_{*}^{-1}\alpha\right)\mu\alpha 
\left(\frac{m+2}{2}\right)^{\frac{2(m+1)}{m+2}}\frac{1}{m+1}
\left(2^{\frac{2(m+1)}{m+2}}-1\right).
\end{align}
Thus, one may rewrite (6.6), (6.7) in a form
\begin{align}
\nonumber
& b_{k}\le C_{b}\left(\frac{\varepsilon_{k}}
{F\Bigl(\lambda \alpha^{-1} h_{*}\xi_{k}^{\frac{m}{m+2}}\Bigr)}\right)^{1/2},\quad
\Delta_{k}\le C_{\Delta}(1-\gamma_{k})\frac{\varepsilon_{k}}
{F\Bigl(\mu \alpha^{-1} h_{*}\xi_{k}^{\frac{m}{m+2}}\Bigr)},\\
& \varepsilon_{k} = \xi_{k+1} - \xi_{k},\quad
F(y) = \frac{\sinh^{2}(y)}{y},\quad
C_{b} = 2\lambda^{1/2}R_{\lambda}\alpha h_{*}^{-1}.
\end{align}
Substituting (6.8) into (5.37), (5.39) one obtains 
\begin{equation*}
\Delta_{k} \le \frac{C_{\Delta}}{C_{b}^{2}C_{L}^{2}}p_{k}^{2}(1-\gamma_{k})
\left(\gamma_{k}-\Delta_{k}\right)^{4}\le 
\frac{C_{\Delta}}{C_{b}^{2}C_{L}^{2}}p_{k}^{2}(1-\gamma_{k})\gamma_{k}^{4}
\end{equation*}
and, consequently,
\begin{equation*}
\gamma_{k+1}\ge \sqrt{1-2 p_{k}}\left(\gamma_{k} - \frac{C_{\Delta}}{C_{b}^{2}C_{L}^{2}}p_{k}^{2}(1-\gamma_{k})^{2}\gamma_{k}^{4}\right).
\end{equation*}
This leads to the following estimate
\begin{align}
\nonumber
& \gamma_{k+1} = \gamma_{0} - \sum\limits_{j=0}^{k}\Delta\gamma_{k},\quad \Delta\gamma_{j} = \gamma_{j} - \gamma_{j+1}\\
\nonumber
& \Delta\gamma_{k}\le \left(\left(1-\sqrt{1 - 2 p_{k}}\right) +
\frac{C_{\Delta}}{C_{b}^{2}C_{L}^{2}} \sqrt{1-2 p_{k}} p_{k}^{2}(1-\gamma_{k})^{2}\gamma_{k}^{3}\right)
\gamma_{k}.
\end{align}
Note that the maximum of the function $(1-x)^{2}x^{3}$ on the interval $[0, 1]$ is attained at $x=3/5$ and equals to $2^{2}3^{3}/5^{5}$. Taking this into account  and the fact that $\gamma_{k+1} \le \gamma_{k}$ for any $k\in \mathbb{N}$, a sufficient condition for convergence of the sequence $\gamma_{k}$ to a positive value $\gamma_{\infty}$ can be written as
\begin{equation}
\sum\limits_{j=0}^{\infty}\left[
\frac{2p_{k}}{1+\sqrt{1-2p_{k}}} + C_{\gamma}\sqrt{1-2p_{k}}p_{k}^{2}\right] < 1,\quad
C_{\gamma} = \frac{C_{\Delta}}{C_{b}^{2}C_{L}^{2}} \frac{2^{2}3^{3}}{5^{5}}.
\end{equation}
One may remark that (6.9) involves only the sequence $\{p_{k}\}_{k=1}^{\infty}$. On the other hand, the conditions (5.37), (5.39), (6.8) imply
\begin{equation}
\frac{\varepsilon_{k}}{F\Bigl(\lambda \alpha^{-1} h_{*}\xi_{k}^{\frac{m}{m+2}}\Bigr)} = 
\frac{p_{k}^{2}\gamma_{k+1}^{4}}{C_{b}^{2}C_{L}^{2}(1-2p_{k})^{2}}
\end{equation}
which defines $\varepsilon_{k}$.

Summarizing all the results and taking into account that $\xi_{k} = \xi_{0} + \sum\limits_{j=0}^{k-1}\varepsilon_{j}$, we arrive at the following theorem
\begin{theorems}
Let $\xi_{0}$ be a positive real number and $\varGamma$ be a non-degenerate geodesic connecting $x_{-}$ and $x_{+}$ on the interval $\Omega_{0}$. It is assumed that $\varGamma$ satisfies the estimate (5.44) and $\xi_{0}$ is small enough such that the estimates (5.43) and (5.47) hold. Assume $h_{*}>0$ satisfies the conditions of Lemma 11 and there exists a sequence $\{p_{k}\}_{k=0}^{\infty}$ such that
\newline
1.\quad $0 < p_{k} < 1/2$ for all  $k$
\newline
2.\quad the condition (6.9) holds
\newline
3.\quad the sequence $\{\xi_{k}\}_{k=0}^{\infty}$ defined by (6.10) has an infinite limit.
\newline
Then there exists a sequence of solutions $q_{k}$ connecting $x_{-}$ and $x_{+}$ on the interval $\Omega_{k}$ and this sequence converges to a transversal doubly asymptotic solution $q_{\infty}$. Moreover, the second derivative of the functional $I$ satisfies
$$
I''[q_{\infty}](\varphi, \varphi) \ge \gamma_{\infty}\Vert \varphi\Vert^{2},\quad
\forall\,\, \varphi\in T_{q_{\infty}}\mathfrak{M}
$$
with $\gamma_{\infty} = \lim\limits_{k\to \infty}\gamma_{k}$.
\end{theorems}

{\bf Remark} We note that a series whose terms are defined by the r.h.s. of (6.10) converges due to assumption (6.9). Besides the function $F$ is strictly increasing on $\mathbb{R}_{+}$. Hence, the condition $\lim\limits_{k\to \infty}\xi_{k} = \infty$ implies
\begin{equation*}
\sum\limits_{k=0}^{\infty}\frac{p_{k}^{2}\gamma_{k+1}^{4}}{C_{b}^{2}C_{L}^{2}(1-2p_{k})^{2}} = 
\sum\limits_{k=0}^{\infty}\frac{\varepsilon_{k}}{F\Bigl(\lambda \alpha^{-1} h_{*}\xi_{k}^{\frac{m}{m+2}}\Bigr)} > 
\int\limits_{\xi_{0}}^{\infty}\frac{{\rm d}\xi}
{F\Bigl(\lambda \alpha^{-1} h_{*}\xi^{\frac{m}{m+2}}\Bigr)}
\end{equation*}
what can be considered as a necessary condition.

\section*{Acknowledgements}

The research was supported by RFBR grant (project No. 17-01-00668/19).

\bigskip

% BELOW SHOULD BE YOUR ADDRESSES
%\address{Alexey V. Ivanov, Saint-Petersburg State University, 199034, Russian Federation, Saint-Petersburg, %Universitetskaya nab. 7/9}

\end{document}